\documentclass[11pt,reqno]{amsart}

\usepackage{amsmath, amsfonts, amssymb, amscd, enumerate, 
}
\usepackage{graphicx}
\usepackage{cyrillic}

\newcommand{\cyrsl}{\fontencoding{OT2}\selectfont\textcyrsl}

\theoremstyle{plain}
\newtheorem{theo}{Theorem}[section]
\newtheorem{lem}[theo]{Lemma}
\newtheorem{prop}[theo]{Proposition}

\newtheorem{cor}[theo]{Corollary}
\theoremstyle{definition}
\newtheorem{rem}[theo]{Remark}

\newtheorem{definition}[theo]{Definition}

\newenvironment{pf}{\noindent{\it Proof. }}{$\square$\par\medskip}

\newenvironment{pfnqed}{\noindent{\it Proof. }}{\par\medskip}

\theoremstyle{plain}

\theoremstyle{definition}

\renewcommand{\=}{\overset{\operatorname{def}}{=}}

\newcommand{\beq}{\begin{equation}}
\newcommand{\eeq}{\end{equation}}
\renewcommand{\a}{\alpha}
\renewcommand{\b}{\beta}

\newcommand{\g}{\gamma}
\newcommand{\h}{\eta}

\renewcommand{\k}{\kappa}
\renewcommand{\l}{\lambda}

\renewcommand{\t}{\tau}

\newcommand{\D}{\Delta}

\newcommand{\G}{\Gamma}

\renewcommand{\L}{\Lambda}
\renewcommand{\O}{\Omega}

\newcommand{\bR}{\mathbb{R}}

\newcommand{\bP}{\mathbb{P}}

\newcommand{\gp}{\mathfrak{p}}

\newcommand{\cC}{\mathcal{C}}

\newcommand{\cF}{\mathcal{F}}
\newcommand{\cG}{\mathcal{G}}

\newcommand{\cP}{\mathcal{P}}

\newcommand{\cS}{\mathcal{S}}

\newcommand{\cU}{\mathcal{U}}

\renewcommand{\square}{\kern1pt\vbox
{\hrule height 0.6pt\hbox{\vrule width 0.6pt\hskip 3pt
\vbox{\vskip 6pt}\hskip 3pt\vrule width 0.6pt}\hrule height0.6pt}\kern1pt}

\renewcommand{\L}{{\cyrsl{L}}} 

\newcommand{\wt}{\widetilde}
\newcommand{\wh}{\widehat}
\newcommand\fro[1]{\stackrel{\displaystyle\frown}{#1}}

\newcommand{\Tang}{\operatorname{\bf t}}
\newcommand{\Norm}{\operatorname{\bf n}}

\newcommand{\n}{\nabla}

\newcommand{\be}{\begin{equation}}
\newcommand{\ee}{\end{equation}}

\def\<#1,#2>{\langle\,#1,\,#2\,\rangle}
\newcommand{\arr}{\begin{array}{rlll}}
\newcommand{\ea}{\end{array}}
\newcommand{\bea}{\begin{eqnarray}}
\newcommand{\eea}{\end{eqnarray}}
\newcommand{\bean}{\begin{eqnarray*}}
\newcommand{\eean}{\end{eqnarray*}}

\def\sideremark#1{\ifvmode\leavevmode\fi\vadjust{
\vbox to0pt{\hbox to 0pt{\hskip\hsize\hskip1em
\vbox{\hsize3cm\tiny\raggedright\pretolerance10000
\noindent #1\hfill}\hss}\vbox to8pt{\vfil}\vss}}}

\newcounter{ssig}
\newcounter{ttig}
\newcommand{\grad}{\operatorname{grad}}

\title[Steepest descent curves on surfaces of constant cvurvature]
{Steepest descent curves   of convex functions \\
 on surfaces of constant curvature}

\author{C. Giannotti and A. Spiro}

   \address
{\newline Cristina Giannotti and Andrea Spiro, 
Dip. Matematica e Informatica, Via Madonna delle Carceri,
I- 62032 Camerino (Macerata),
ITALY\newline
\phantom{a}}

\email
{cristina.giannotti@unicam.it}\par
\medskip

 \email
{andrea.spiro@unicam.it}\par
\medskip
\date{\today}


\begin{document}

\begin{abstract} Let $\cS$ be a complete surface of constant curvature $K = \pm 1$,  i.e.  $S^2$  or $\L^2$, and $\O \subset \cS$  a bounded convex subset.  If $\cS = S^2$, assume also $\operatorname{diameter}(\O) < \frac{\pi}{2}$.  It is proved that the  length of  any  steepest descent curve of  a quasi-convex function in  $\O$ is  less than  or equal to the perimeter of $\O$. This upper bound is actually proved  for the class of $\cG$-curves, a family of curves that naturally includes all steepest descent curves. In case $\cS = S^2$,  it is also  proved  the existence of $\cG$-curves, whose  length is equal to the perimeter of their convex hull, showing that  the above estimate is indeed optimal. The results generalize theorems by Manselli and Pucci on steepest descent curves in the Euclidean plane.  
\end{abstract}
\maketitle
\section{Introduction}
\setcounter{equation}{0}
Let $\cS$ be a complete surface of constant Gaussian curvature $K = 0$, $+1$ or $-1$,  i.e.  the Euclidean plane $E^2$, the unit sphere $S^2$ or the  Lobachevskij plane  $\L^2$.  
An absolutely continuous curve $\g: [a, b] \to \cS$ (e.g.  a rectifiable curve parameterized by arc length)  is called   {\it $\cG$-curve}  if it  satisfies the following condition
 (here and in the rest of the paper,  we use the notation $\g_s = \g(s)$): \par
\smallskip
\moveright-0.7cm
\vbox{
\begin{itemize}
\item[]{\it  for any $s_o \in [a,b]$  for which 
$\dot \g_{s_o}$ exists and is different from $0$, 
 all points $\g_{s}$ with $s \leq s_o$   are  in the same  closed half space,   bounded    by  the normal   to $\g$ at $\g_{s_o}$\/}.
 \end{itemize}}
 \smallskip
 \noindent  \par
\noindent Notice that the  class of   $\cG$-curves  naturally includes all  steepest descent curves of convex functions, that is the  $\cC^1$-curves $\g_t$ satisfying  equations of the form 
 $\dot \g_t =  - \grad f|_{\g_t}$  for some convex  $\cC^1$-function  $f: \cU \subset \cS \to\bR$. 
The $\cG$-curves have been originally considered   by Manselli and Pucci in \cite{MP},  where they determined an optimal  upper bound for the length of $\cG$-curves
contained in  a given   bounded  convex subset of $E^2$. \par
\medskip
 Here  we consider the problem of  establishing similar   upper bounds for the lengths of $\cG$-curves in $S^2$ and $\L^2$. Our first result  consists of  the following generalization of the upper bound determined in  \cite{MP}: \par
\begin{theo} \label{first} Let  $\cS = S^2, E^2$ or $\L^2$ and $\g: [a,b] \to \cS$ a $\cG$-curve. In case  $\cS = S^2$,
assume also  that  $\operatorname{diameter}( \g([a,b])) < \frac{\pi}{2}$. Then the length  $\ell(\g)$ of $\g$ is less than or equal to the perimeter  $\gp(\wh\g)$ of the convex hull $\wh \g$ of $\g$. 
\end{theo}
Notice that, as in Euclidean geometry, also  in $S^2$ or $\L^2$ the perimeter of a bounded convex set is less than  or equal to the perimeter of any larger convex set  (Prop. \ref{Proposition1}). Due to this,  all $\cG$-curves in a given convex set have length less than  or equal to the perimeter of that set.  \par
\medskip
We also remark that the condition on  $\operatorname{diameter}( \g([a,b])) $ in case $\cS = S^2$ is a technical hypothesis, needed in our analysis of  the  growth of the perimeters $\gp(s)$ of 
the convex hulls $\wh{\g|_{[a,s]}}$  of the arcs $\g|_{[a,s]}$, $s \in (a,b]$. It can be replaced by other assumptions, applicable when $\operatorname{diameter}( \g([a,b])) \geq \frac{\pi}{2}$  (like e.g. that  the tangent to $\g$ at $\g_s$ is a tangent line also to the convex hull $\wh{\g|_{[a,s]}}$ for any $s$). But,    at the moment,  the previous statement is the best and most general one, which  we have at our reach. \par
\medskip
Secondly, for checking the optimality of this estimate, we consider the existence problem for  $\cG$-curves  $\g: [a, b] \to \cS$ {\it with maximal length property},  i.e. such that
the length    $\ell(\g|_{[a,s]})$   equals  to the perimeter of the convex hull $\wh {\g([a,s])}$  for any $s \in [a,b]$.  The solution of this problem in  case  $\cS = E^2$  was the second main result in 
   \cite{MP} (see also \cite{MP1}),  where  the authors proved that, up to rigid motion,  for given $L$ there exists    a unique $\cG$-curve,  of  length $L$ and $\cC^1$ outside  the starting point, with  the  maximal length property. This curve  is the unique curve in $E^2$ which  is {\it self-involute\/}, i.e.  it is equal   to the initial arc  of  its own involute  (see \S  4 for definitions; see also e.g. \cite{Gu, Ra}). \par
\smallskip
From the  proof of  Theorem \ref{first}, it is natural to conjecture  that  a similar result should be true also  when   $\cS = S^2$ or $\L^2$, at least for  sufficiently small curves.  Our second main result  shows that for the  case $\cS = S^2$  this expectation is indeed correct.  In fact,  we first prove  that  {\it there exists on $S^2$ curves which are self-involute}Ê (Corollary \ref{thecorollary})  and then we prove that {\it  any self-involute of $S^2$ is also a  $\cG$-curve with maximal length property},  showing in this way the  optimality of the upper bound of Theorem \ref{first}. \par 
\medskip
The existence problem for  $\cG$-curves with maximal length property on $\L^2$ remains open. It might  be also interesting to know if the self-involutes of $S^2$ share the same uniqueness properties of the self-involutes of $E^2$ and whether there exists or not self-involutes on $\L^2$. 
\par
\medskip
The structure of the paper is as follows.  In \S 2, we collect a few basic facts on convex sets in surfaces of constant curvature. For reader's convenience, we tried to make the exposition as much as possible self-contained. In \S 3, we introduce the notion of $\cG$-curves    and prove Theorem \ref{first}. In \S 4, we recall some well-known  facts of the Differential Geometry of  curves on surfaces, we introduce the definitions of involutes, almost self-involutes and self-involutes and prove our final results on the existence of self-involutes and $\cG$-curves with maximal length property on $S^2$.\par
\medskip
\noindent{\it Acknowledgment.} The authors   would like to thank Paolo Manselli, to whom are very much indebted  for suggesting  the problem, giving  most of the   ideas   for the proof of Theorem \ref{existenceselfinvolutes} and  for a multitude of very cheerful and  fruitful discussions on this subject. We are also grateful to Gabriele Bianchi and the referee for valuable  comments and helpful remarks. \par
\bigskip
\section{Preliminaries}
\setcounter{equation}{0}
\label{prelim}
\subsection{Convex sets in surfaces of constant curvature}\hfill\par
\medskip
\label{preliminaries1}
Consider  an (abstract) complete  surface    $\cS$  of constant curvature $K
= +1, 0$ or  $-1$,  i.e. a complete 2-dimensional Riemannian manifold locally
isometric to  the  unit sphere $S^2$, the Euclidean
plane $E^2$ or the  Lobachevskij plane  $\L^2$.  It is well-known that either
 $S^2$ or $E^2$ or $\L^2$ is  also the universal
covering space of $\cS$.\par
\smallskip
A  {\it line of $\cS$}  is  the trace of a
maximal geodesic and,  for any two points $x_1$, $x_2 \in \cS$, 
the   {\it (geodesic) segment joining  $x_1$ and $x_2$}  is an arc of
minimal length between those points.  We can now recall the definition of
convex sets in $\cS$ (see e.g. \cite{Al}). 
 \smallskip
 \begin{definition}\label{convexitydef}
A  subset $Q \subset \cS$ is called {\it convex} if for any two points $x, x' \in Q$ 
there exists a unique  segment joining $x$ to  $x'$ and such segment is  
entirely  in $Q$.  
Given a subset $A \subset \cS$,  its {\it convex hull} is the intersection  $\wh A = \bigcap A'$ of all
convex sets   $A' \subset \cS$ containing $A$.\end{definition}
 Notice  that, according to this definition,  there are  sets with no convex hull. For instance, if $A \subset \cS$ contains a pair of  points with  more than one segment  joining them,  there is no convex set that contains $A$.\par
\smallskip
We also recall that, according to  how many segments  might join  two given points and whether or not  they are all  included in the subset,  other notions of convex subsets in a surface can be given (see e.g. \cite{Al, BZ, Ud}).  However  all these notions coincide when $\cS = E^2$ or $\L^2$, while for the  case  $\cS = S^2$  the above choice turned out to be the most convenient  one for establishing an upper bound for the length of  steepest descent curves. 
\par
\medskip
 Let $\pi: \wt \cS \longrightarrow \cS$
be the universal covering space of  $\cS$. 
If $\cU \subset \cS$ is an open
convex set, then it is  simply connected and the restriction  $\pi|_{\wt \cU}: \wt
\cU  \longrightarrow \cU$ of $\pi$ to a connected component $\wt \cU$ of $\pi^{-1}(\cU)$ is an isometry between  the
 convex subset $\wt \cU$ of $\wt \cS$ and  $\cU$.  Hence,  for our purposes, with no loss of generality  we
may always reduce  to the cases   $\cS = S^2, E^2$ or $\L^2$.    In addition, 
the following lemma shows that {\it  the case  $\cS = S^2$  can be always replaced
by the assumption 
$$\cS = S^2_+ = S^2 \cap \{\ (x^1,x^2,x^3)\ :\
x^3
> 0\ \}\ .$$}
\par 
\begin{lem} \label{lemma22}ÊLet $\cU \subset S^2$ be open and convex. Then it is contained in a
hemisphere.
\end{lem}
\begin{pf} Intuitively,  the claim is a consequence of the fact that a convex subset of the unit sphere cannot contain a pair of antipodal points.  But a short and precise proof can be obtained as follows.  Let $C\subset \bR^3$ be the cone given by the half-lines from the
origin passing through the points of $\cU$.  One can check that $\cU$ is
convex if and only if $C$ is convex. If we denote by $\a \subset \bR^3$ a
supporting plane of $C$ through the origin, it follows immediately that 
$\cU = S^2 \cap C$ is contained in a hemisphere bounded by $\a$. 
\end{pf} 
Consider the upper hemisphere $S^2_+$. 
The map $\varphi: S^2_+ \to  \{\ x^3= 1\ \} \simeq \bR^2$, 
sending  any   $x \in S^2_+$ into its  radial projection $\bar x \in  \{\ x^3=
1\ \} $,   maps the  lines of $S^2_+$ into the lines of the Euclidean space $E^2 = (\bR^2, g_o)$ 
(here $g_o$ is the standard Euclidean metric). 
In fact,  the traces of geodesics
in $S^2_+$ are   great circles,  i.e.   intersections between $S^2_+$
and  affine planes of $\bR^3$ 
 through the origin,  and are  mapped by $\varphi$  into 
 straight lines, given by the intersections of those  planes with $\{\ x^3=
1\ \}$.\par  \medskip
Similarly, consider the Lobachevskij plane 
$$\L^2 = (\{\ (x^1)^2 + (x^2)^2 -  (x^3)^2 = - 1\ ,
\ x^3> 0\ \}, g)\ ,$$ 
where $g$ is the Riemannian metric induced from the
indefinite metric on $\bR^3$   
$$g = dx^1 \otimes dx^1  + dx^2 \otimes dx^2 -  dx^3 \otimes dx^3\ . $$
Also in this case,  the lines of $\L^2$ are given by the intersections
between $\L^2$ and the affine planes of $\bR^3$ 
 through the origin.  So, the map  $\psi: \L^2 \longrightarrow  \{\ x^3= 1\
\} \simeq \bR^2$,   sending the points of   $\L^2$ into their radial 
projections from the origin on that plane, maps $\L^2$ into  the open unit disc
$\D \subset \bR^2 $  and the lines of $\L^2$ into the 
straight lines (chords)  of $\D$.  \par
We recall that  the abstract surface $ K^2 = (\Delta,
\wt g \= \psi^{-1}{}^*(g))$  is usually called {\it Klein disc}.\par
\bigskip
Using the maps $\varphi$ and $\psi$,  we may always
identify $S^2_+$ and $\L^2$ with (an open subset of) $\bR^2$, 
 endowed with a suitable metric $g = g_{ij} dx^i \otimes dx^j$, whose
geodesic segments  coincide with  the standard  Euclidean segments.  Under this
identification,  {\it a subset  $Q \subseteq \cS$ is convex
 if and only if it is convex   in the Euclidean sense.}
\par
For this reason, in the  figures of this paper, the 
segments and  convex subsets of  $\cS$  will be drawn as   Euclidean segments
and Euclidean convex subsets of $\bR^2$. However, the reader should be aware
that,  when $\cS$ is not the Euclidean space,  these picture can  be misleading for what concerns  lengths and
angles. For instance,  for any $v, w \in T_x \cS \simeq T_x E^2$, the norm $|v|$ and the cosine $\cos(\wh{v w})$ are in general different from the
Euclidean values, since they are given   in terms of the non-Euclidean metric 
$g = g_{ij} dx^i \otimes dx^j$ of $\cS$ by   $$|v| = \sqrt{g_{ij}(x) v^i v^j } \
,\qquad \cos(\wh{v w}) = \frac{g_{ij}(x) v^i w^j}{|v| |w|}\ .$$ 
\par
\medskip
A closed simple curve  $\cP \subset \cS$ is called {\it polygonal\/} 
if  it is  union of finitely many segments,  called {\it sides\/}. 
The  endpoints of the sides are called {\it vertices\/}. \par 
The distance  between  two points $x_o, y_o \in \cS (\subseteq \bR^2)$  is
equal to  the length of the segment between   $x_o$ and $y_o$ w.r.t.  $g$ i.e. $$d_\cS(x_o, y_o) = \int_0^1\sqrt{
g_{ij}(x(t)) \dot x^i(t) \dot x^j(t)} dt \ ,$$  where $x(t)$ is the line $$x(t) = (x_o^1 (1-t) + t y_o^1, x_o^2 (1-t) + t y_o^2)\ .$$  
In particular, if $x_o, y_o$ are in  a fixed convex compact subset $K \subset 
\cS$,  there are $0 < C_K, C'_K$ such that $C_K |x_o - y_o|  \leq d_\cS(x_o, y_o)
\leq C'_K |x_o - y_o|$.  \par
\smallskip
The 
length $\ell_\cS(\cP)$ of a polygonal curve $\cP$ is  the sum of  the lengths
of  its sides. A curve $\cC \subset \cS$ is called {\it rectifiable} if its 
length $$\ell_\cS(\cC) = \sup\{\ \ell_\cS(\cP)\ ,\ \cP\ \text{polygonal curve with
vertices in}Ê\  \cC\ \} $$ is finite. For an absolutely continuous parameterization $\g:
[a,b] \to \cC \subset  \cS $ of a rectifiable curve $\cC$,  the tangent vector
$\dot \g_{s_o}$ exists for almost any $s_o \in[a,b]$ and
$$\left.\frac{d \ell_\cS(\g|_{[a,s]})}{ds}\right|_{s_o} = |\dot \g_{s_o}| =
\sqrt{g(\dot \g_{s_o}, \dot \g_{s_o})}\ .$$
By previous remarks, for any compact subset $K \subset \cS (\subset \bR^2)$,  
  a curve $\cC \subset K$ is rectifiable if and only if it is rectifiable 
as curve in the Euclidean plane and 
\begin{equation} \label{estimate1} C_K\cdot  \ell(\cC) \leq \ell_\cS(\cC) \leq C_K' \cdot \ell(\cC) \end{equation}
where $\ell(\cC)$ is the length of $\cC$ in Euclidean sense.\par
\bigskip
The following proposition  generalizes two well-known properties of convex sets in $E^2$.\par
\begin{prop} \label{Proposition1} \hfill\par
\begin{itemize}
\item[i)]  The boundary $\partial Q$ of a  bounded convex set   $Q \subset \cS$ is  a rectifiable curve. 
\item[ii)] If   $Q \subset
Q'$ are two bounded convex subsets of $\cS$, then $\ell_\cS(\partial Q) \leq
\ell_{\cS}(\partial Q')$. \end{itemize} \end{prop}
\begin{pf} (i) Since   $Q \subset \cS \subset \bR^2$ is  convex
also in the  Euclidean sense, the claim follows immediately from   known facts
on Euclidean convex sets and \eqref{estimate1}.\par \smallskip
 (ii) 
Let $\cP$ be a polygonal
curve with vertices $x_1, \dots, x_n, x_{n+1}Ê= x_1 \in \partial Q$ and 
$\bP$ the  (open) convex polygon with $\partial \bP = \cP$.  
Denote also by  $s_i$ the  side of $\cP$  joining $x_{i}$ and $x_{i+1}$.  \par
 Being $Q'$  bounded and convex, the line $\wh s_1$
containing $s_1$  cuts  $Q'$ into  two convex subsets and the polygon $\bP$ is entirely contained  in  one of them (see Fig.÷1).  Call   $Q'_{(1)}$ this convex set and observe that  $\partial Q'_{(1)}$ is obtained by replacing  a portion of  $\partial Q'$ with  the segment joining  the endpoints  of such portion.  Hence  $\ell_{\cS}(\partial Q') \geq \ell_{\cS}(\partial Q'_{(1)})$. Next,   consider the line $\wh s_2$ containing the side $s_2$ and the  convex subset  $Q'_{(2)} \subset Q'_{(1)}$,  which is cut by $\wh s_2$ and contains $\bP$. As before, we have that  $\ell_{\cS}(\partial Q'_{(1)}) \geq \ell_{\cS}(\partial Q'_{(2)})$.  
Repeating the same  construction for all the  lines  $\wh s_i$ containing  the sides $s_i$,  we end up   with a nested sequence of convex sets $Q'_{(i)}$, $1 \leq i \leq n$,   with $Q'_{(n)} = \bP$  and such that 
$$\ell_{\cS}(\partial Q') \geq \ell_{\cS}(\partial Q'_{(1)}) \geq \ell_{\cS}(\partial Q'_{(2)})  \dots \geq  \ell_{\cS}(\partial \bP) = \ell_\cS(\cP)$$ 
By arbitrariness of   $\cP$,  it follows that  $\ell_{\cS}(\partial Q') \geq \ell_{\cS}(\partial Q)$. 
\end{pf}
\leftline{\includegraphics[width=3.5cm]{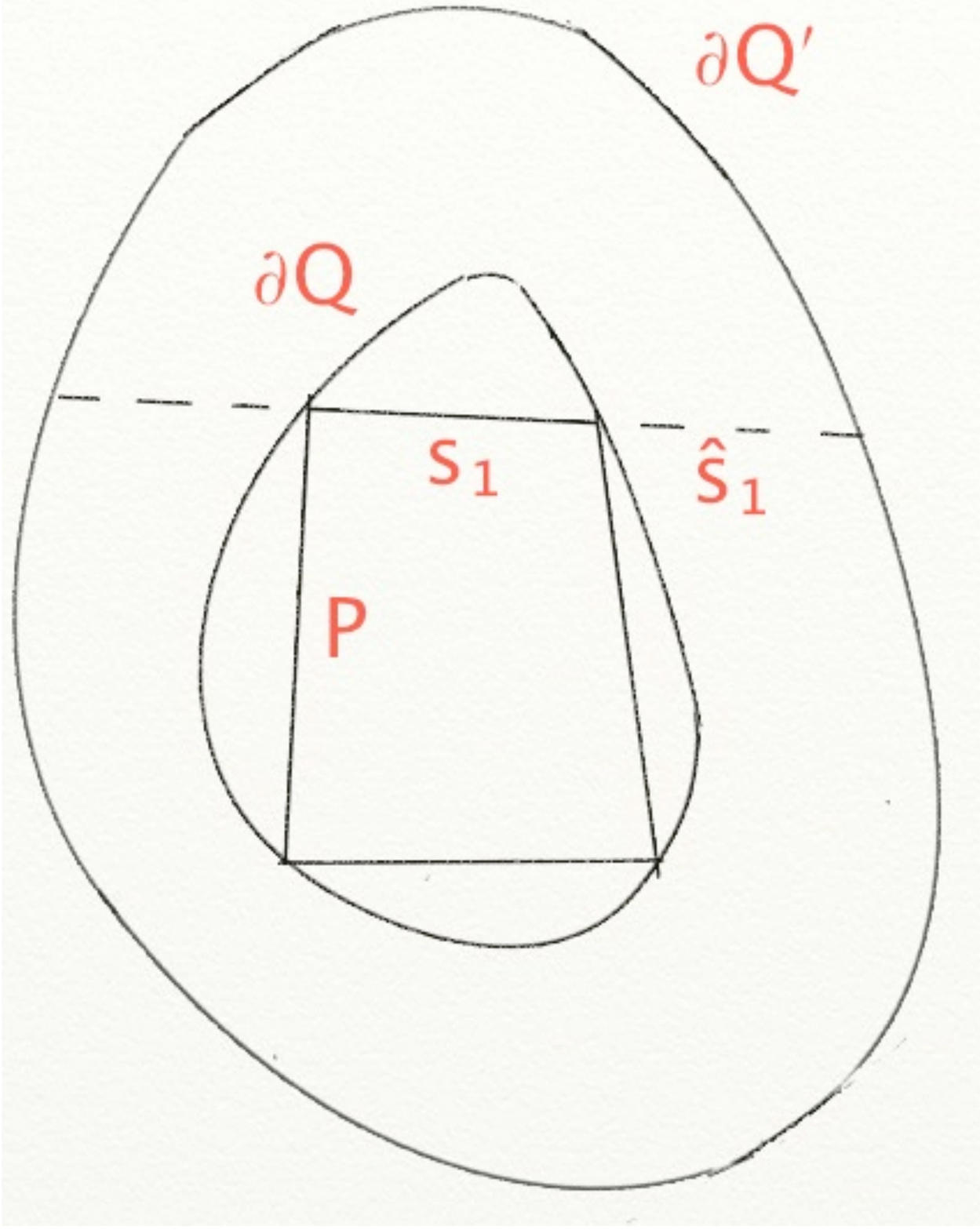}}
\vskip -5 truecm
\centerline{\includegraphics[width=4cm]{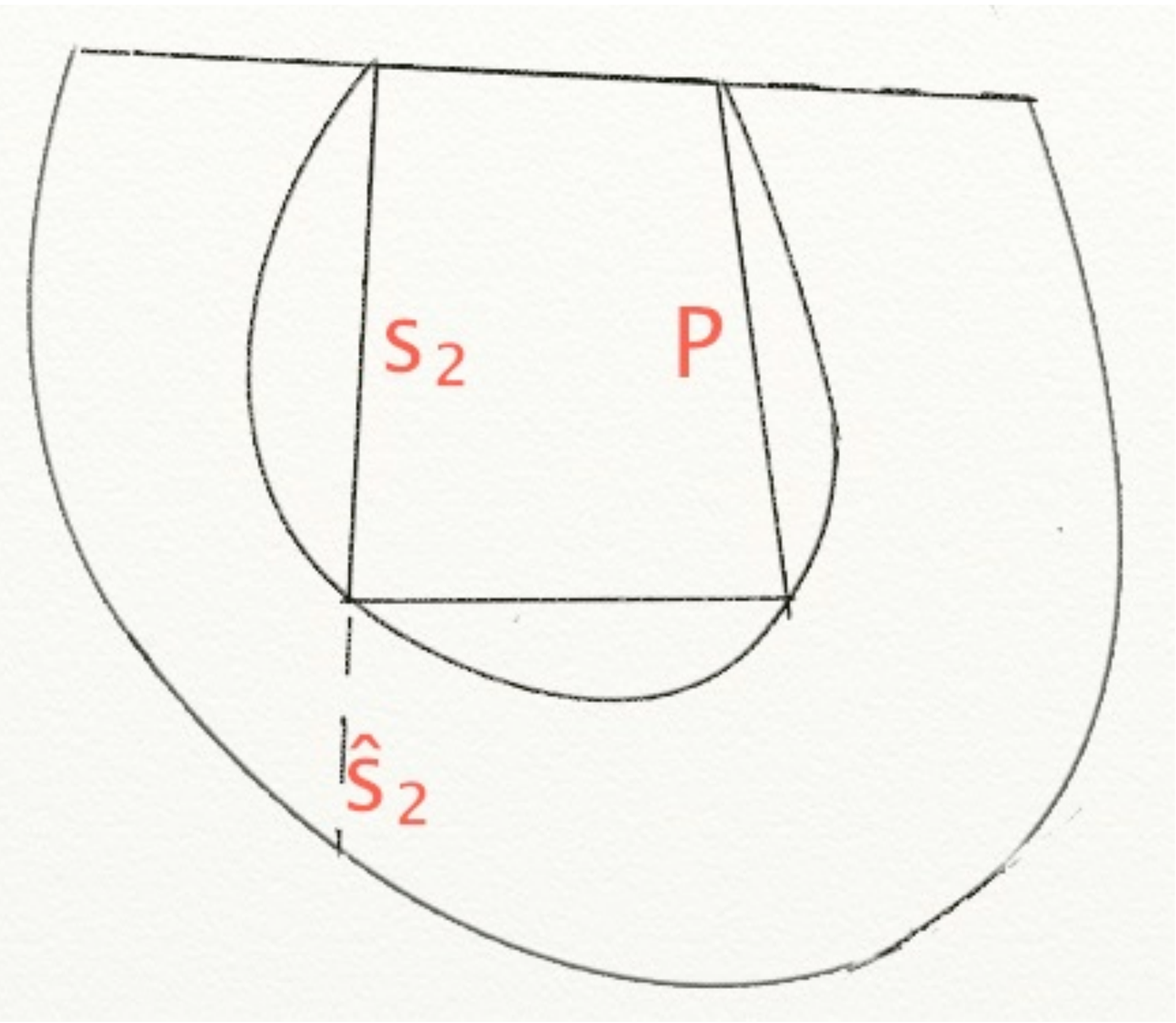}}
\vskip -5 truecm
\rightline{\includegraphics[width=3cm]{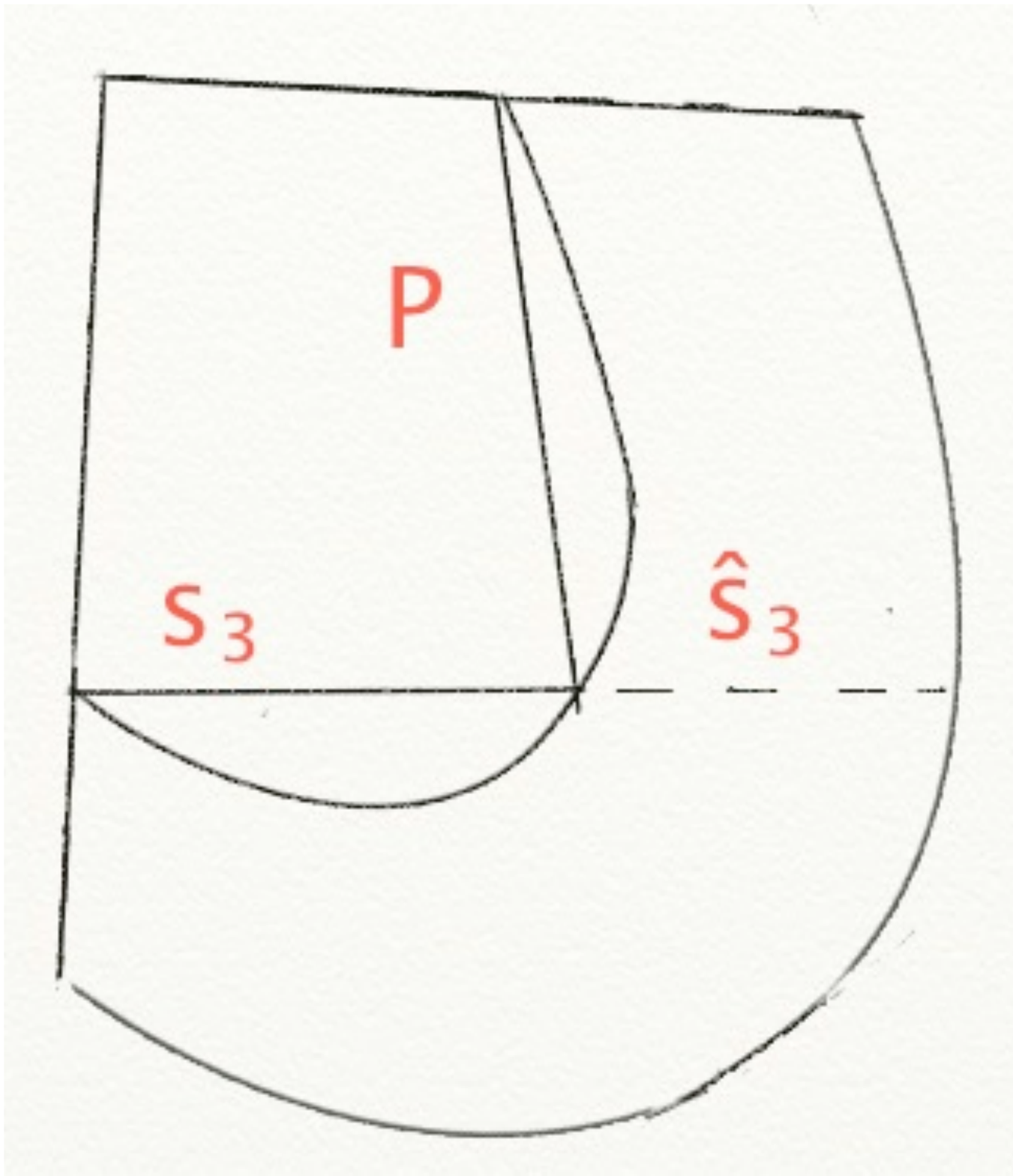}}
\vskip - 0.5 truecm
\centerline{\Small \bf Fig. 1}
\begin{rem} By a refinement of  the proof, it is not hard to check  that if $Q, Q'$ are as in   Proposition \ref{Proposition1} (ii) and if  $Q' \setminus   Q$ has non empty interior, then $\ell_\cS(\partial Q') > \ell_{\cS}(\partial Q)$.  
\end{rem}
Given a bounded convex set $Q \subset \cS$,  we  denote
by  $\gp(Q)$ the  perimeter $\ell_\cS(\partial Q)$.\par
\medskip 
\subsection{An  auxiliary lemma}
\hfill\par
\medskip
We give here a lemma,  needed in  the proof of Theorem 1.1. It could be  derived from  the Gauss Lemma on derivatives of distance functions and  it  could  be proved  for  Riemannian manifolds. However, we provide a   self-contained proof, using only basic facts on surfaces of  constant curvature.\par
\medskip
In the following,  given  $y, z \in \cS = S^2_+$, $E^2$ or
$\L^2$, we denote by $[y,z]$
the segment joining them and, given  
 a curve $\eta: [a,b] \to \cS$, we denote by
$$d_{y,\eta}: [a,b] \to \bR\ ,\qquad d_{y,\eta}(s) = \ell_\cS([y,\eta_s]) =
d_\cS(y, \eta_s)
\ .$$ 
\begin{lem} \label{lemmino} Let $\eta: [a,b] \to \cS$ be a curve in $\cS$,   $s
\in [a,b[$,  such that  the tangent vector $\dot \eta_s \neq 0$ exists at $\eta_s$,  and 
$\{h_n\}$  a sequence of  positive (negative) real numbers converging to $0$ and
$\{y_n\} \subset \cS$ a sequence of points, all of them different from $\h_s$ and   $\h_{s + h_n}$,  converging to   a point $y_o$.  
\par
We denote  by $\a_n$ and $\g_n$ the angles, with vertices in $\h_s$ and $y_n$, respectively, 
formed by the tangent vectors  of the oriented segments $[\h_s, \eta_{s+ h_n}]$ and $[\h_s, y_n]$  and by the tangent vectors of the oriented segments 
$[y_n, \eta_{s+ h_n}]$ and $[y_n, \h_s]$.  
If  $y_o \neq \h_s$ 
\begin{equation} \label{formulalemmino} \lim_{n \to \infty} \frac{d_{y_n,
\eta}(s + h_n) - d_{y_n, \eta}(s)}{h_n} = |\dot \eta_s| \cos \varphi\ .  \end{equation}
where $\varphi =  \lim_{n \to \infty}Ê(\pi - \a_n)$. This equality holds also  when  $y_o = \h_s$, provided that 
$\lim_{n \to \infty} \alpha_n$ exists and  $\lim_{n \to \infty} \g_n = 0$. \par
Moreover,  $d_{y_o, \eta}'(s) = |\dot \eta_s| \cos \varphi$ for any $y_o \in \cS$, where 
$\varphi$ denotes  the angle   between $\dot \h_s$ and  the tangent in $\h_s$  of the oriented segment $[\h_s, y_o]$ when $y_o \neq \h_s$ and 
$\varphi = 0$ when $y_o = \eta_s$.
\end{lem}
\begin{pf} First of all, let us recall  the following well-known formulae of
Spherical and Hyperbolic Trigonometry (see e.g. \cite{AVS} \S I.3).  \par
\smallskip
\noindent{\it On the sphere $S^2$:}  
$$
  \cos \a  = - \cos \b \cos \g +
\sin \b \sin \g \cos a\ ,$$
 $$ \frac{\sin \a}{\sin a} = \frac{\sin
\b}{\sin b} = \frac{\sin \g}{\sin c}\ ,$$
\par
\smallskip
\noindent{\it On the Lobachevskij plane $\L^2$:} 
$$\cos \a = - \cos \b \cos \g + \sin \b
\sin \g \cosh a \ ,$$
$$ \frac{\sin \a}{\sinh a} = \frac{\sin
\b}{\sinh b} = \frac{\sin \g}{\sinh c}\ ,$$
where, $a$, $b$, $c$ are
the sides of a triangle and $\a$, $\b$ and $\g$ the corresponding opposite 
angles.  Now, consider the case  $h_n > 0$ for all $n$ (the case $h_n < 0$ is similar).   For a given $n$, consider the triangle of vertices $y_n$, $\eta_s$
and $\eta_{s+h_n}$ and let us denote by
$$a_n = \ell_\cS(
[y_n, \eta_{s+h_n}])\ ,\qquad b_n = \ell_\cS( [y_n , \eta_{s}])\
,\qquad c_n = \ell_\cS([\eta_s, \eta_{s+h_n}])$$
and  by $\a_n$, $\b_n$ and $\g_n$ the corresponding opposite  angles. 
 Then,   both the hypotheses  $  y_o  \neq \h_s$ and   $y_o = \h_s$  with  $\lim_{n \to \infty} \a_n = \pi - \varphi$, $\lim_{n \to \infty} \g_n = 0$,   imply that
 \begin{equation}\label{formula3bis} \lim_{n\to\infty}\left(a_n - b_n
\right)=0\ ,\qquad   \lim_{n\to\infty} c_n=0\ ,
\end{equation}
\begin{equation} \label{formula3ter} 
\lim_{n\to\infty}\g_n=0\ ,\quad
\lim_{n\to\infty}\a_n = \pi - \varphi\ ,\quad 
\lim_{n\to\infty}\b_n=\varphi\ .\end{equation}
Observe that 
$$\lim_{n \to \infty} \frac{d_{y_n, \eta}(s + h_n) - d_{y_n, \eta}(s)}{h_n} =
\lim_{n \to \infty}Ê\frac{a_n - b_n}{c_n} \cdot \lim_{n \to \infty}
\frac{c_n}{h_n}\ .
$$
The second limit is equal to $|\dot \eta_s|$, while the first limit can be written as 
\begin{equation}\label{form4} \lim_{n\to \infty} \frac{a_n - b_n}{c_n} =
\lim_{n\to \infty} \frac{\sin(a_n - b_n)}{\sin(c_n)} = \lim_{n \to \infty} \frac{\sin a_n \cos b_n - \cos a_n
\sin b_n}{\sin c_n}\end{equation}
 or as 
\begin{equation}\label{form4bis} \lim_{n\to \infty} \frac{a_n - b_n}{c_n} =
\lim_{n\to \infty} \frac{\sinh(a_n - b_n)}{\sinh(c_n)} = \lim_{n \to \infty}
\frac{\sinh a_n \cosh b_n - \cosh a_n \sinh b_n}{\sinh c_n}\end{equation}
Using  \eqref{formula3bis},  \eqref{formula3ter} and the above relations in 
\eqref{form4}  when $\cS = S^2_+$    or in \eqref{form4bis} when $\cS =
\L^2$,  we  obtain  that in both cases
 $$ \lim_{n\to \infty} \frac{a_n -
b_n}{c_n} =   \lim_{n \to \infty}Ê \frac{\left(\cos \b_n - \cos
\a_n\right)\left(1 - \cos  \g_n\right) }{\sin^2
\g_n}  = \cos \varphi\ .$$
In case $\cS = E^2$, using the Sine Law of Euclidean trigonometry and the
 standard property of Euclidean triangles $\a_n + \b_n + \g_n = \pi$, 
 we have  $$ \lim_{n\to \infty} \frac{a_n -
b_n}{c_n} =  \lim_{n\to \infty} \frac{\sin \a_n -
\sin\b_n  }{\sin \g_n} = $$
$$ = \lim_{n\to \infty} \frac{ 2 \sin\left(\frac{\a_n - \b_n}{2} \right)
\cos\left(\frac{\a_n + \b_n}{2} \right)}{\sin
\g_n} = \lim_{n\to \infty} \frac{ 2 \sin\left(\frac{\a_n - \b_n}{2} \right)
\sin\left(\frac{\g_n }{2} \right)}{\sin
\g_n} = \cos \varphi\ ,$$
where the last equality follows from \eqref{formula3ter}.\par
For what concerns  the last claim, in case $y_o \neq \h_s$ it is a direct consequence of  \eqref{formulalemmino} for the sequence $\{y_n = y_o\}$, while in case $y_o =\h_s$ it follows from  the immediate observation  that $d'_{y_o, \h}(s) = \lim_{h\to 0} \frac{d_\cS(\h_s, \h_{s+h})}{h} = | \dot \h_s|$. 
\end{pf}
\section{$\cG$-curves}
\setcounter{equation}{0}
\bigskip

Let  $f: \cU \subset \cS \longrightarrow \bR$ be a map defined on a convex subset $\cU$ of $\cS$. 
We say that $f$ is  {\it quasi-convex} if any   level set 
$\cU_{\leq c} \= \{ \ x \in \cU\ : \ f(x) \leq c\ \}$ is   convex in $\cS$.
\footnote{{\it Warning}. In \cite{Ud} and other places,  ``quasi-convex function''  
means a  function on  a totally convex set  with totally convex level  sets. 
Notice  that if $\cS = E^2$,  the notion of convexity and total convexity coincide.}
Notice that, under the identification of $\cS$ with $\bR^2$ or $\D \subset
\bR^2$, $f$ is quasi-convex as function on $\cS$  if and only if it is
quasi-convex  in the usual Euclidean sense.\par \medskip

Let  $f: \cU \subset \cS\longrightarrow  \bR$ be  a $\cC^1$ 
quasi-convex function $f$ with $df \neq 0$ at all points, so that
 the  level sets $\cU_{\leq t}$ have $\cC^1$-boundaries.  Any  
curve $\g_s$ of {\it steepest descent\/} for $f$  (i.e.   such that   $\dot \g_s =  -\left.\grad f\right|_{\g_s}$) 
intersects orthogonally all the boundaries  of the level sets and,  for given
$s_o$,   all points $\g_{s}$, $s \geq  s_o$, are included  in the level set $ \cU_{\leq f(\g_{s_o})}$. 
 This means that the class of  steepest descent curves  is naturally included
in the following class of curves, which extends the class   considered in \cite{MP}.\par 
\begin{definition} We say that an absolutely continuous  curve $\g: [a, b] 
\to \cS$ is  {\it in the class $\cG$\/} if, for any $s_o \in [a,b]$ such that
$\dot \g_{s_o}$ exists and is different from $0$, 
 all points $\g_{s}$ with $s \leq s_o$ are in  a same closed half space that
is
bounded    by  the normal  line $\ell$  to $\g$ at $\g_{s_o}$.\par
  \end{definition}
 \medskip 
In the following, we  denote by  $\g: [a,b] \to \cS$ a $\cG$-curve of $\cS$ 
and, for any $s\in [a,b]$ we  indicate by $\gp(s)$ the perimeter of the convex
hull of $\g|_{[a,s]}$. By Proposition \ref{Proposition1} (ii), the function $\gp(s)$ is
not decreasing and hence $\gp'(s)$ exists for almost all $s \in [a,b]$.\par
Moreover, for any   $x \in \cS$ and $v, w \in T_x \cS$, we will denote by
$C(x; v,w)$ a closed convex  sector bounded by the two geodesic rays originating from
$x$ and tangent to $v$ and $w$. Identifying   $\cS$ with (an open subset of) $\bR^2$, the sector $C(x; v,w)$  is  the convex angle with vertex $x$ and sides parallel to $v$ and $w$. It is clearly always uniquely determined  
except when  it is a half-plane, i.e. when  $\wh{v w} = \pi$.  \par
\smallskip
For any point $\g_s$ of the curve $\g$, 
we call {\it projecting sector of $\g$ at $\g_s$} the smallest closed convex sector
   containing $\g|_{[a,s]}$. For any point $\g_s$ of $\g$, we will denote by 
   $v_i = v_i(s) \in T_{\g_s} \cS$, $i = 1,2$,  the tangent  vectors of the boundary rays of  the corresponding  projecting  sector, which will be therefore indicated by $C(\g_s; v_1, v_2)$. \par
\begin{theo} \label{maintheorem1} Let  $\g: [a,b] \to \cS$ be 
a  $\cG$-curve  with convex hull $\wh \g$ and, in the case   $\cS = S^2_+$,
assume also  that $\operatorname{diameter}(\g) < \frac{\pi}{2}$.
Then 
\beq\label{fundamentalinequality}  \ell_\cS(\g) \leq
\gp(\wh\g) \ ,\eeq
where the equality  occurs only if for almost  all $s \in [a,b]$ the
projecting sector $C(\g_s; v_1, v_2)$ of $\g$ is such that 
$\wh{v_1 v_2} = \frac{\pi}{2}$ and  either $v_1$ or $v_2$ is tangent to $\g$ at
$\g_s$. \par
In particular, for any $\cG$-curve  in a bounded convex set $Q \subset \cS$
(and satisfying   the above condition  when $\cS = S^2_+$) the length 
$\ell_\cS(\g)$ is less than  or equal to $\gp(Q)$.  \end{theo} 
The proof of this result is based on the following two lemmata.\par
\begin{lem} \label{lemma1}ÊLet $s$ be such that both $\gp'(s)$ and  $\dot \g_{s}$
  exist with $\dot \g_s \neq 0$. Then
$$\frac{\gp'(s)}{|\dot\g_{s}|} \geq \cos \phi_1 + \cos \phi_2\ ,
 \qquad \text{with}\ \ \phi_i \= \pi -\wh{\dot\g_{s} v_i}\ ,$$ 
where $v_1$, $v_2$ are the  vectors of  the projecting sector
$C(\g_{s}; v_1, v_2)$.
\end{lem}
\begin{pfnqed} For  $s \in [a,b[$ and $h > 0$, consider the following
notations:     
\begin{itemize}
\item[--] $A_s = \g|_{[a,s]}$ and  $\wh{A_s}$ is  its convex
hull; 
\item[--]  $A_s(h) = A_s \cup \{\g_{s + h}\}$ with convex
hull  $\wh{A_s(h)}$; 
\item[--]  $\fro {A_s (h)} = C(\g_{s}; v_1, v_2) \cap \wh{A_s(h)}$ .
\end{itemize}
\par
\centerline{\includegraphics[width=8cm]{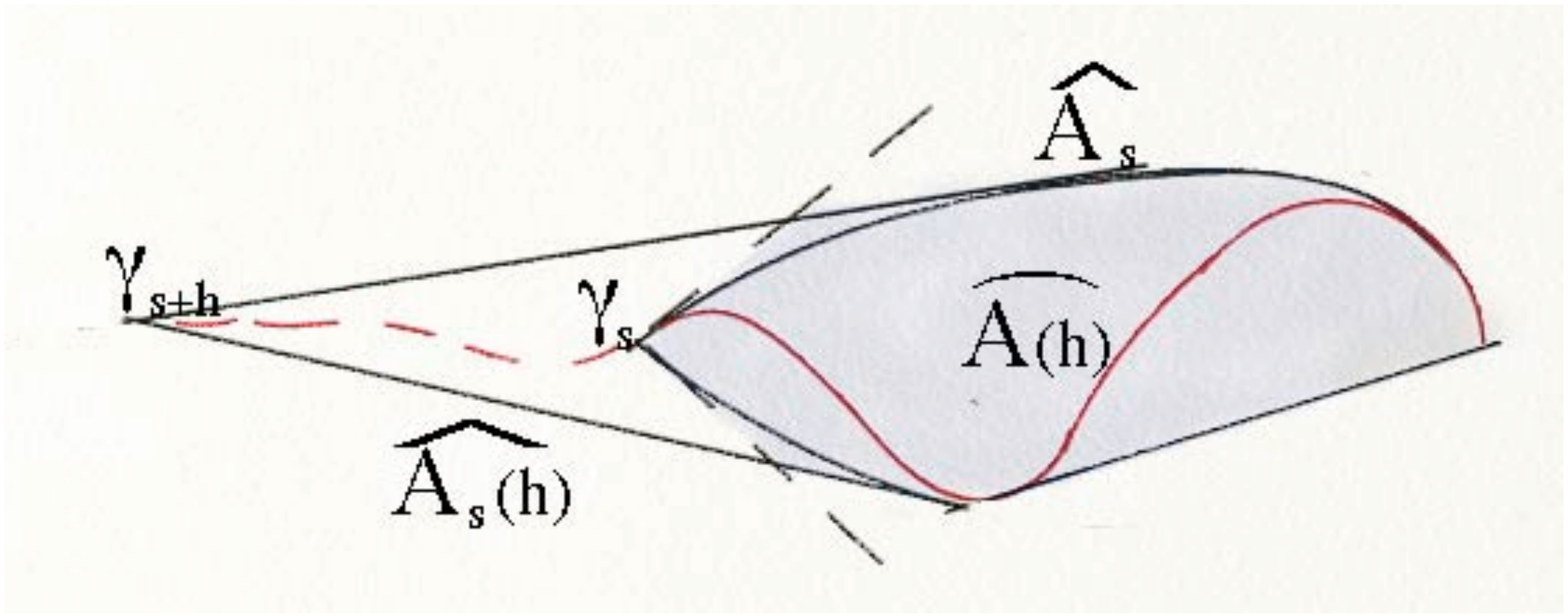}}
\vskip-1 truecm
\centerline{\Small \bf  Fig. 2}
Clearly, 
$$\wh{A_s}  \subseteq \fro{A_s(h)} \subseteq \wh{A_s(h)} \subseteq 
\wh{A_{s+h}} $$
and hence, by Proposition \ref{Proposition1} (ii), $\gp(\wh{A_s}) \leq 
\gp(\fro{A_s(h)})$ $\leq \gp(\wh{A_s(h)}) \leq
\gp(\wh{A_{s+h}})$  and
\begin{equation}\label{form1}\gp(s+h) - \gp(s) = \gp(\wh{A_{s+h}}) -
\gp(\wh{A_s}) \geq \gp(\wh{A_s(h)}) - \gp(\fro{A_s(h)})\ . \end{equation} 
The boundary of $\wh{A_s(h)}$ contains two segments,  lying on two rays coming
out from $\g_{s+h}$, and these segments necessarily intersect the sides of
$C(\g_s; v_1, v_2)$ in two distinct points, which we call $x_1^h$ and
$x_2^h$ (one of them might be $\g_s$). Moreover, since the set of points $\{x_i^h\ ,\ h>0\}$ is bounded, there
exists a sequence $\{h_n\}$ with $\lim_{n\to\infty}h_n=0$  and
$\lim_{n\to\infty}x_i^{h_n}=x_i$  for some  $x_i$   in
one  side of $C(\g_{s}; v_1, v_2)$.\par
Using the notation of Lemma \ref{lemmino}, we may write
$$ \gp(\wh{A_s(h_n)}) - \gp(\fro{A_s(h_n)}) =  \sum_{i =
1}^2\left(\ell_\cS( [x_i^{h_n}, \g_{s+h_n}]) - \ell_\cS( [x_i^{h_n} , \g_{s}])
\right) = $$
\begin{equation} \label{form2} = \sum_{i = 1}^2 \left(d_{x_i^{h_n}, \g}(s+ h_n) -
d_{x_i^{h_n}, \g}(s)\right) \ . \end{equation}
Recall that, for fixed $i = 1,2$, 
 all $x_i^{h_n}$  lie  in one of the two sides (call them $r_i$, $i = 1,2$)  of $C(\g_s; v_1, v_2)$. 
Hence, if  there exists  a  subsequence   $\{x_i^{h_{n_k}}\}$ converging to  a point $y_o  \neq \g_s$,  then $y_o \in r_i$   and 
the angles $\a_{n_k}$,  formed by the tangents in $\g_s$ of the oriented segments $[\g_s, \g_{s + h_{n_k}}]$ and $[\g_s,  x^{h_{n_k}}_i]$,   converge to $\pi - \phi_i$. By  Lemma \ref{lemmino} we get (after replacing $\{x_i^{h_n}\}$ by  the above  subsequence)  
\beq \label{f1} \lim_{n \to \infty}\frac{d_{x_i^{h_{n}}, \g}(s+ h_{n}) - d_{x_i^{h_{n}},
\g}(s)}{h_{n}} = |\dot \g_s| \cos \phi_i\ .\eeq 
The same conclusion holds also if there exists  a  subsequence   $\{x_i^{h_{n_k}}\}$ of points, all different from $\g_s$,   converging to   $y_o  = \g_s$, because in such case, by possibly taking another subsequence,  the angles in $x^{h_{n_k}}_i$, formed by  the oriented segments $[x_i^{h_{n_k}}, \g_{s + h_{n_k}}]$ and $[x_i^{h_{n_k}},\g_s]$, tend to $0$ and  Lemma \ref{lemmino} applies.  \par
 To check this claim, consider the rays $r^{n_k}_i$ with origin in $\g_{s + h_{n_k}}$ and  containing $[\g_{s+ h_{n_k}} ,  x^{h_{n_k}}_i]$. By construction, any such ray lies  in a line of support  for $\g|_{[a, s]}$ and contains a sub-ray,  with origin  in $x^{h_{n_k}}_i$,  included in 
$ C(\g_s; v_1, v_2)$. By taking a suitable subsequence,  the rays  $r^{n_k}_i$  tend to a ray $r^o_i$, with origin in $\g_{s}$, 
entirely included in $C(\g_s; v_1, v_2)$ and lying  in a line of support $\ell_o$ for $\g|_{[a, s]}$. 
Due to this,  if  $r^o_i$ intersected the interior of  $C(\g_s; v_1, v_2)$,   this    would not be the smallest  convex sector containing $\g|_{[a, s]}$. 
So $r^o_i \subset \partial C(\g_s; v_1, v_2)$, that is  $r^o_i = r_i$,  and  the angles in $x^{h_{n_k}}_i$, formed by   $[x_i^{h_{n_k}}, \g_{s + h_{n_k}}]$ and $[x_i^{h_{n_k}},\g_s]$,  tend to $0$ as claimed. \par
\smallskip
Now, we claim that \eqref{f1} is true also if  there is no  subsequence $\{x_i^{h_{n_k}}\}$, made   of points all different from $\g_s$. In fact, in this case  we may assume that $x_i^{h_n} = \g_s$ for any $n$ and hence that  the projecting sector $C(\g_s; v_1, v_2)$  lies in the intersection of   half spaces, bounded by the lines $\ell_n$, which contain the boundary segment $[\g_{s + h_n}, \g_s] \subset \partial \wh{A_s(h_n)}$.  We remark that 
\begin{itemize}
\item[a)]  the lines $\ell_n$ tend to the tangent line $\ell_o$ of $\g$ at  $\g_s$ and 
$C(\g_s; v_1, v_2)$ is contained in a half-space bounded by $\ell_o$; 
\item[b)]Êsince  the rays, with origin in  $\g_s$ and containing $[\g_s, \g_{s - h}]$, $h > 0$, are included in  $C(\g_s; v_1, v_2)$ and tend to a ray  $r_o \subset \ell_o$ when $h \to 0$, it  follows that  $r_o \subset C(\g_s; v_1, v_2)$. 
\end{itemize}
From  this,  we infer that  $r_o = r_i\subset  \partial C(\g_s; v_1, v_2) \cap \ell_o$ and $\phi_i =  \pi - \wh {\dot \g_s v_i}  = 0$.  So, by Lemma \ref{lemmino},
$$\lim_{n \to \infty}\frac{d_{x_i^{h_{n}}, \g}(s+ h_{n}) - d_{x_i^{h_{n}},
\g}(s)}{h_{n}} = d'_{\g_s, \g}(s) = |\dot \g_s| \cos  0 = |\dot \g_s| \cos  \phi_i$$ 
as claimed. From \eqref{form1} and \eqref{form2}, 
$$\gp'(s) \geq \sum_{i = 1}^2 
\lim_{n \to \infty}\frac{d_{x_i^{h_n}, \g}(s+ h_n) - d_{x_i^{h_n},
\g}(s)}{h_n} = |\dot \g_s| (\cos \phi_1 + \cos \phi_2)\ .\eqno\qed$$
\end{pfnqed}
\medskip
\begin{lem} Let $s$ be such that  $\dot \g_{s}$ exists with $\dot \g_s \neq 0$ and let $C(\g_{s}; v_1, v_2)$
be the projecting sector of $\g$ at $\g_{s}$. 
If the curvature of $\cS$ is $K = + 1$, assume also that
$\operatorname{diameter}(\g|_{[a,s]}) < \frac{\pi}{2}$.
 Then
$\wh{v_1 v_2} \leq \frac{\pi}{2}$  and hence, if  $\phi_i = \pi - 
\wh{\dot\g_{s} v_i}$, 
$$\cos \phi_1 +
\cos \phi_2 \geq 1\ .$$ 
The equality holds if and only if $\wh{v_1 v_2} = \frac{\pi}{2}$ and one of the
$\phi_i$'s is equal to $0$. \end{lem}
\begin{pf} First of all, we claim that, for any $s_o \in [a,b[$, the function
$$d_{\g_{s_o}, \g}:[s_o, b]\to \bR \ ,\qquad d_{\g_{s_o}, \g}(s) = d
_\cS(\g_{s_o}, \g_s)$$ is non-decreasing. In fact,  since $\g$ is a $\cG$-curve, the vector
$\dot \g_s$ and the tangent vector in $\g_s$ to the oriented 
segment $[\g_{s_o}, \g_{s}]$ points towards  the same  side w.r.t. the
normal line of $\g$ in $\g_s$. In particular, the angle $\varphi$ between them
is less than or equal to $\frac{\pi}{2}$. By Lemma \ref{lemmino}, 
$d_{\g_{s_o}, \g}'(s)$  is non-negative for almost all $s$ and the
claim follows. \par
\medskip
Secondly, we claim that for any $a \leq s_1 < s_2 < s \leq b$, the  
angle $\a$ formed by the oriented segments $[\g_{s}, \g_{s_1}]$ and
$[\g_{s}, \g_{s_2}]$ is less than  or equal to $\frac{\pi}{2}$. \par
In case $\cS =
E^2$ or $\L^2$,  it can be checked as follows.  In the triangle with vertices
$\g_{s}$, $\g_{s_1}$ and $\g_{s_2}$, the sum of inner angles is less than or equal to $\pi$ 
and hence $\a > \frac{\pi}2$   only if it is the largest of these
three angles.  But this cannot be because the  side 
$[\g_{s_1}, \g_{s_2}]$,  opposite to $\a$, is not the largest one (by the previous claim, 
it  is shorter or equal to 
$[\g_{s_1}, \g_{s}]$), 
in contrast with a well known fact of Euclidean and Hyperbolic Geometry.\par
Also in case the curvature of $\cS$ is $K = +1$ and $\operatorname{diameter}(\g|_{[a,s]}) <
\frac{\pi}{2}$,  if $\a$ were larger than $\frac{\pi}{2}$, its opposite side in the triangle with vertices
$\g_{s}$, $\g_{s_1}$ and $\g_{s_2}$ would be the largest one, as it  can be checked  using the spherical law of cosines $\cos a = \cos b \cos c + \sin a \sin b \cos \a$. Hence, also in this case we conclude that $\a \leq \frac{\pi}{2}$ by the same argument as above. \par
Now, the first statement of the lemma  follows immediately 
from   the observation  that  $\wh{v_1 v_2} = \phi_1 + \phi_2$   is limit
 of angles delimited by    segments  $[\g_{s}, \g_{s_1}]$ and
$[\g_{s}, \g_{s_2}]$ for some $a \leq s_1 < s_2 < s \leq b$. To check the
second statement, just look for    the   minimum  of   $f(\phi_1,
\phi_2) = \cos \phi_1 + \cos \phi_2$ in the region  $\O = \{\ 0 \leq \phi_i \leq
\frac{\pi}{2}\ , \phi_1 + \phi_2  \leq \frac{\pi}{2}\ \} $.
 \end{pf}
Combining the results of these lemmata,  $\gp'(s) \geq |\dot \g_s|$ for almost
all $s \in [a,b]$,   with equality only  if  the
amplitude of the  projecting sector is $\frac{\pi}{2}$, and this implies  the theorem.\par
 \bigskip
 
 \section{Self-Involutes  on spheres}
 \setcounter{equation}{0}
 \subsection{Basic facts on curves of surfaces of constant curvature}\hfill\par
In this section $\cS$ is a simply connected, complete surface of constant curvature with a curvature $K$ that might   assume any real  value.  Recall that,   when $K = \pm 1/ R^2$, the surface $\cS$  is either
$$S^2_R = \{\ x  \in \bR^3 \ : \ x^T \cdot  x  = R^2\ \} $$
or 
$$ \L^2_R =\!\! \left\{\ x \in \bR^3 \ : \  x^3 > 0 \ \text{and} \ \ x^T  \cdot I_{2,1} \cdot x  = - R^2\ \text{where} \ I_{2,1}Ê= 
   \left(\smallmatrix 1 & 0 & 0 \\ 0 & 1 & 0 \\
 0 & 0 & -1
 \endsmallmatrix\right)\ \right\},Ê$$
 respectively endowed with the Riemannian  metric $g$, which is  induced  either by the standard Euclidean metric  $g^E = dx^1 \otimes dx^1  + dx^2 \otimes dx^2 +  dx^3 \otimes dx^3$  or  by the Lorentzian metric  $g^L = dx^1 \otimes dx^1  + dx^2 \otimes dx^2 -  dx^3 \otimes dx^3$. For uniformity of notation,  one can consider also  $E^2$  as a surface in $\bR^3$, namely setting $E^2 = \{\ x  \in \bR^3 \ :\ x^3 = 1 \ \}Ê\subset \bR^3$, with the metric induced  by the standard Euclidean metric $g^E$ of $\bR^3$. \par
 \medskip
 Given a $\cC^3$ curve $\eta: I \subset \bR \to \cS \subset \bR^3$, parameterized by arc length,  the {\it Frenet frame of $\h$\/} at $s = s_o$ is the orthonormal basis  $(\Tang_s, \Norm_s) \subset T_{\h_s} \cS$,  given by $\Tang_s = \dot \h_s$ and the unit vector $\Norm_s$, tangent to $\cS$, orthogonal to $\Tang_s$ and so that $(\h_s, \Tang_s, \Norm_s)$ is a positively oriented basis for $\bR^3$.  By  a simple generalization of the classical theory of  curves in $E^2$,  one can derive the following ``Frenet formulae''  
\beq \label{Frenetformulae} \n_{\dot \h_s} \Tang_s = \kappa_s \Norm_s \ ,\qquad \n_{\dot \h_s} \Norm_s = - \kappa_s \Tang_s\ ,\eeq
 for some smooth function $\k_s$, called {\it (geodesic) curvature\/} of $\h$. We also recall the Levi-Civita connection $\n$ of $\cS$ is such that, for any smooth curve $\h$ and any vector field $Y_{\h_s}$, tangent to $\cS$ and  defined at the points of $\h$,
 $$\n_{\dot \h_s} Y = \dot Y_{\h_s} +  K g(\dot  \eta_s, Y_{\h_s}) \h_s$$
 (here,  $ \dot Y_{\h_s}$ and  plus sign denote the standard first order derivative and the sum of  maps   with values in $\bR^3$). Also, given a point $x_o \in \cS \subset \bR^3$ and a unit vector $v \in T_{x_o} \cS$, the geodesic $\g_s$ with $\g_0 = x_o$ and $\dot \g_s = v$ is of the form
$$\g_s = \cos\left(\frac{s}{R}\right) x_o + \sin\left(\frac{s}{R}\right) (R v) \quad \text{or}\quad  \g_s = \cosh\left(\frac{s}{R}\right) x_o + \sinh\left(\frac{s}{R}\right) (R v)$$
\beq  \text{or}\qquad \g_s = x_o + s v\eeq
 according to the value of  $K$.\par
 \smallskip
 We conclude this subsection, stating the   following  simple generalization  of  the classical  ``Fundamental Theorems of Plane Curves''.  It is an immediate consequence of   the Existence and Uniqueness Theorem  for O.D.E.'s.\par 
\begin{theo}\label{fundamentaltheorem}\hfill\par
\begin{itemize}
\item[1)] Let $\h: [0,L] \subset \bR \to \cS$ and $\h': [0,L'] \subset \bR \to \cS$ be two $\cC^2$ curves parameterized by arc length and  with curvature functions  $\k$ and $\k'$, respectively.  There exists an isometry $g: \cS \to \cS$ such that $\h' = g \circ \h$ if and only if $L = L'$ and $\k_s = \k'_s$ for any $s \in [0,L]$. 
\item[2)]  Let  $\k:  [0, L] \to \bR$ be $\cC^0$. For any  $x_o \in \cS$,   $v \in T_{x_o} \cS$ with $|v| = 1$ and  $\wt L \leq L$ sufficiently small,  there exists a unique $\cC^2$ curve $\h: [0, \wt L) \to \cS$ parameterized by arc length and with curvature function $\k$,  such that $\h_0 = x_o$ and $ \dot \h_0 = v$.
\end{itemize}
\end{theo}
\medskip
\subsection{Involutes  and (almost) self-involutes}
\begin{definition}  Let $\h: [0, L] \to \cS \subset \bR^3$ a curve,  parameterized by arc length (hence $\dot \h_s = \Tang_s$) and $\cC^1$ on $(0,L)$.   The {\it involute of $\h$}  is the curve $\wt \h:   [0, L] \to \cS$   such that, for any $s \in (0,L)$, the point $\wt \h_s \in \cS$ is determined by the formula
\beq \wt \h_s = \g_s^{(\h_s, - \Tang_s)}\ ,\eeq
where  $\g^{(\h_s, - \Tang_s)}$ is the geodesic of $\cS$ with $ \g_0^{(\h_s, - \Tang_s)} = \h_s$ and 
$\dot \g_0^{(\h_s, - \Tang_s)} = - \Tang_s$. \par
The curve $\h_s$ is called {\it almost self-involute} if there exists an isometry of $\cS$, which maps $\h$ into the initial arc of  length $L$  of its own involute $\wt \h$ (re-parameterized by arc length). In case $\h_s$ coincides with such initial arc of $\wt \h$, we call it {\it self-involute}.
\end{definition}
\leftline{\includegraphics[width=6.5cm]{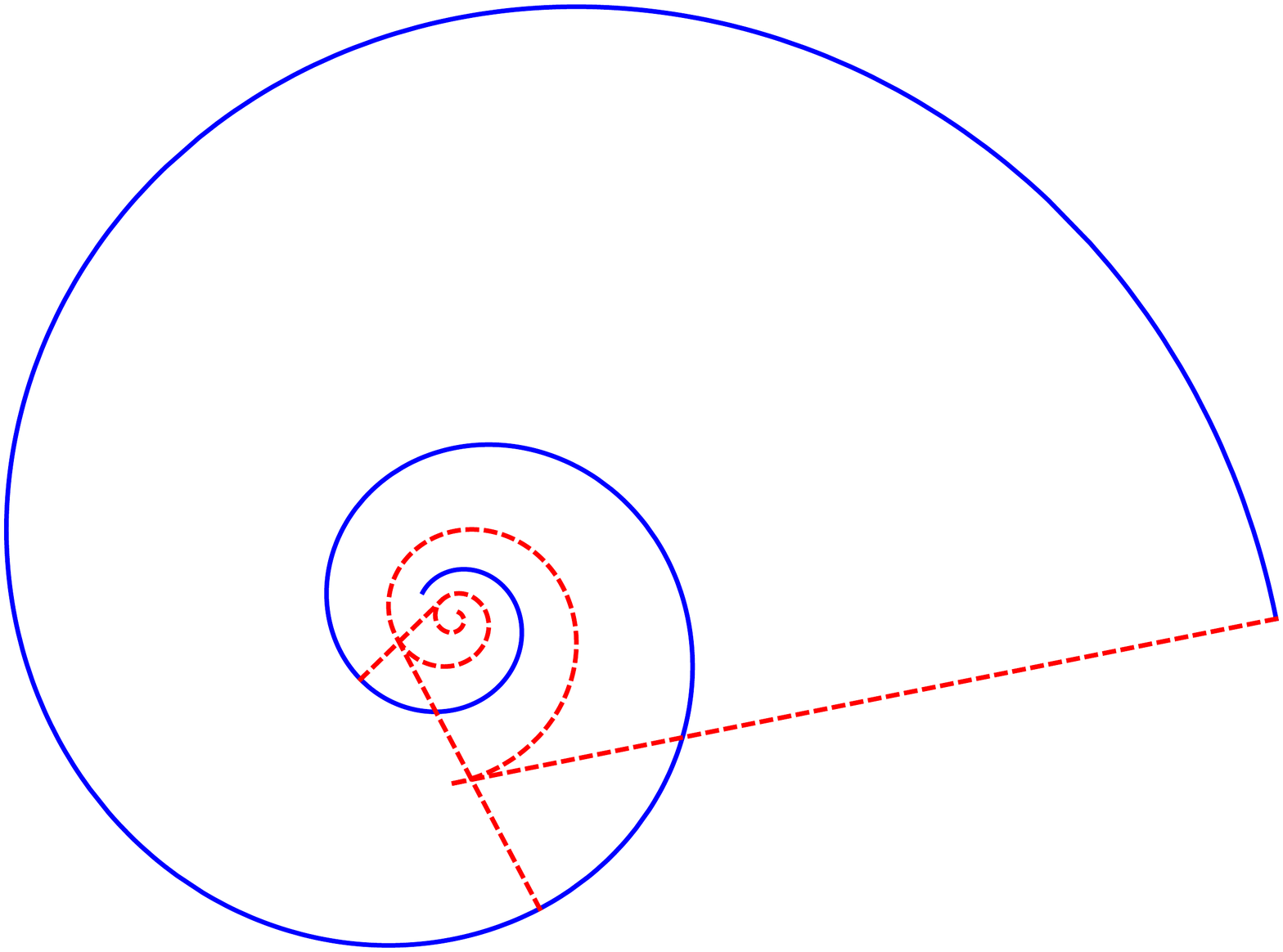}}
\vskip - 4.56 truecm
\rightline{
\includegraphics[width=6.5cm]{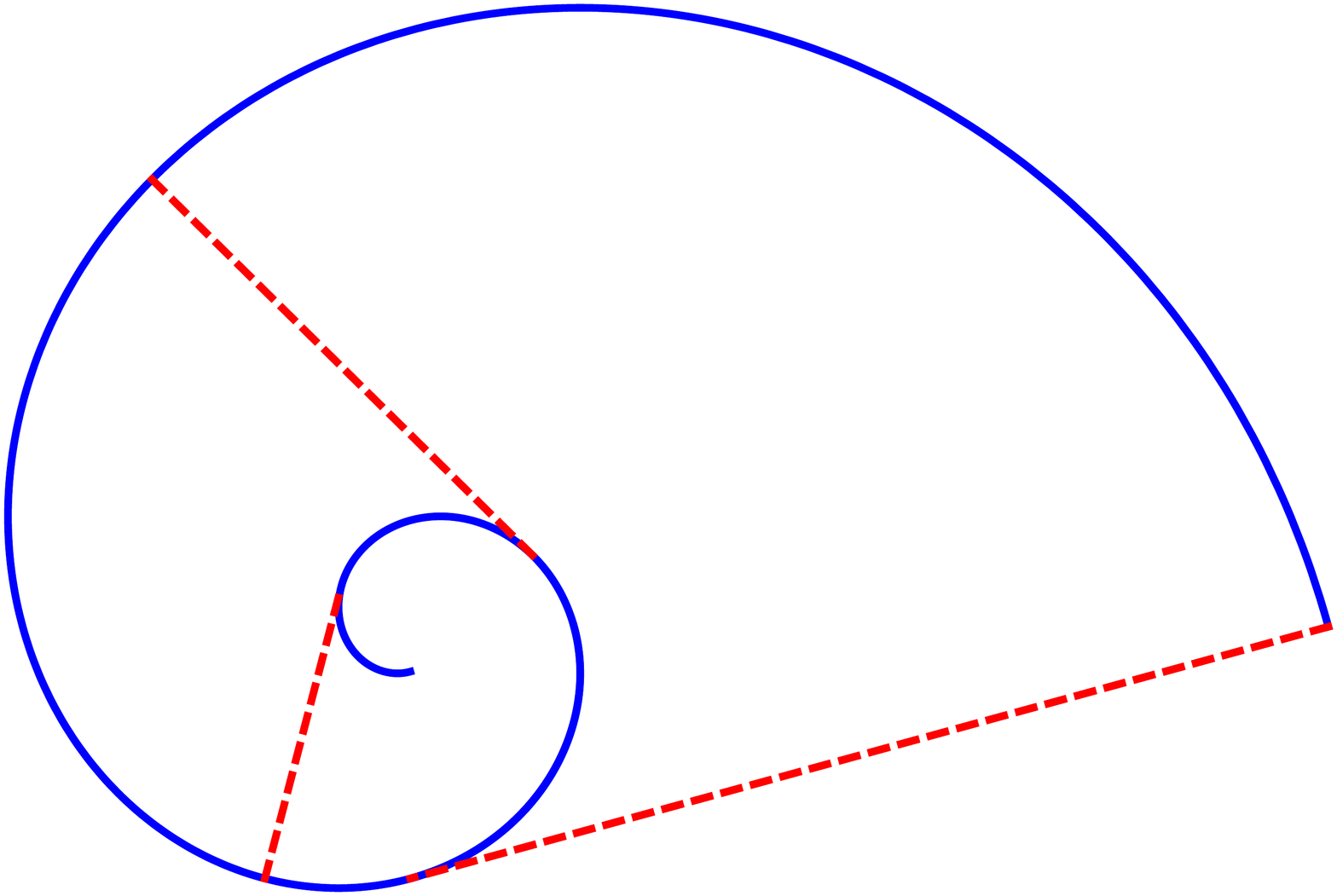}}
\centerline{\Small \bf Fig. 3 -  An almost self-involutes (left) and a self-involute (right)}
\medskip
\par
According to the value of $K$, one immediately obtains that
$$\wt \h_s =  \cos\left(\frac{s}{R}\right) \h_s - \sin\left(\frac{s}{R}\right) (R \dot \h_s) \quad \text{or}\quad  \wt \h_s = \cosh\left(\frac{s}{R}\right) \h_s - \sinh\left(\frac{s}{R}\right) (R \dot \h_s)$$
\beq \label{involute}  \text{or}\qquad \wt \h_s = \h_s -  s \dot \h_s\ .\eeq
Simple computations shows that, in case $\h$ is $\cC^2$,  the Frenet frames $(\Tang_{\wt s_s}, \Norm_{\wt s_s})$, the arc length parameter $\wt s$  and the curvature $\wt k_{\wt s_s}$ at  the points $\wt \h_{\wt s_s} = \wt \h_s$ are (only for $s < \pi/2$ when $K > 0$ and up to  changes $\wt s \to - \wt s$)  
\beq \label{4.14} \Tang_{\wt s_s} = - \Norm_s\ , \quad \Norm_{\wt s_s} = \left\{\begin{array}{ll}  \sin\left(\frac{s}{R}\right) \frac{ \h_s}{R}   + \cos\left(\frac{s}{R}\right) \dot \h_s & \smallmatrix \text{if} \ K =  \frac{1}{R^2} \endsmallmatrix\\
\phantom{a}& \phantom{a}\\    \sinh\left(\frac{s}{R}\right) \frac{ \h_s}{R}  + \cosh\left(\frac{s}{R}\right)  \dot \h_s&  \smallmatrix \text{if} \ K = - \frac{1}{R^2} \endsmallmatrix\\ \phantom{a}& \phantom{a}\\ \dot \h_s  &\smallmatrix \text{if} \ K = 0 \endsmallmatrix\end{array} \right.\ ,
\eeq
$$\dot{\wt s}_s = \left\{\begin{array}{ll} R \sin\left(\frac{s}{R}\right) \k_s & \smallmatrix \text{if} \ K =  \frac{1}{R^2} \endsmallmatrix\\
\phantom{a}& \phantom{a}\\  R \sinh\left(\frac{s}{R}\right) \k_s&  \smallmatrix \text{if} \ K = - \frac{1}{R^2} \endsmallmatrix\\ \phantom{a}& \phantom{a}\\ s  \k_s &\smallmatrix \text{if} \ K = 0 \endsmallmatrix\end{array} \right.\ ,\quad
 \wt \k_{\wt s_s} = \left\{\begin{array}{ll} \frac{1}{R} \cot\left(\frac{s}{R}\right)&\smallmatrix \text{if} \ K =  \frac{1}{R^2} \endsmallmatrix\\
\phantom{a}& \phantom{a}\\  \frac{1}{R} \coth\left(\frac{s}{R}\right) &\smallmatrix \text{if} \ K = - \frac{1}{R^2} \endsmallmatrix \\ \phantom{a}& \phantom{a}\\ \frac{1}{s} &\smallmatrix \text{if} \ K = 0 \endsmallmatrix\end{array}\right.$$
Using Theorem \ref{fundamentaltheorem}, one can infer that a curve $\h:[0,L] \to S^2_R$, parameterized by arc length, smooth and with strictly positive curvature function $\k_s$ on $(0,L)$, is congruent to the initial arc of length $L$ of its own involute  if and only if there exists a strictly increasing smooth function $\t: (0, \wt L) \subset (0,L)Ê\to (0,L)$, with $\lim_{s \to 0^+} \t_s = 0$ and satisfying the following equations
\beq \label{system} \left\{\begin{array}{l} \dot \t_s = R \sin\left(\frac{s}{R}\right) \k_s  \\ \phantom{a}\\ \k_{\t_s} = \frac{1}{R} \cot\left(\frac{s}{R}\right)\ . \end{array}\right.\eeq
Moreover,  for any pair of smooth functions $\k: (0,L) \to \bR$ and $\t: (0, \wt L) \subset (0,L) \to (0,L)$ with $\k_s, \dot \t_s  > 0$ for any $s \in (0,L)$,   $\lim_{s \to 0^+} \t_s = 0$  and satisfying \eqref{system}, there exists a curve $\h:[0,L] \to S^2_R$ with curvature $\k_s$ and  congruent to the initial arc of its own involute.  Similar claims hold for the curves in $\L^2_R$ or $E^2$, where \eqref{system} is  replaced by 
\beq \label{system1} \left\{\begin{array}{l} \dot \t_s = R \sinh\left(\frac{s}{R}\right) \k_s  \\ \phantom{a}\\ \k_{\t_s} = \frac{1}{R} \coth\left(\frac{s}{R}\right) \end{array}\right.\qquad \text{or} \qquad  \left\{\begin{array}{l} \dot \t_s = s \k_s  \\ \phantom{a}\\ \k_{\t_s} = \frac{1}{s}  \end{array}\right.\ ,\eeq
respectively. \par
\medskip
The curvature functions $\k_s$ of  almost self-involutes are exactly  those which  satisfy \eqref{system} or \eqref{system1} with  suitable $\t_s$ and the cardinality (up to congruences) of the self-involutes of a given length is clearly bounded by the cardinality of the solutions to those systems.  \par
\medskip
\subsection{An existence result for the almost self-involutes on spheres}\hfill\par
\begin{theo} \label{existenceselfinvolutes}ÊFor any $a > 1$ and any $s_o > 0$ sufficiently small, there  exists an almost 
self-involute  $\h: [0, s_o] \to S^2_R$ of class $\cC^3$, with curvature $\k_s$ and related function $\t_s$, which satisfy \eqref{system} on  $(0, s_o]$ and  such that 
\beq \label{5.17} \lim_{s\to 0^+} \t_s = 0\ ,\quad  \lim_{s\to 0^+} \dot \t_s = a\ ,\quad  \lim_{s\to 0^+} \k_s = + \infty\ . \eeq
\end{theo} 
\begin{pf} By Theorem \ref{fundamentaltheorem} and the remarks in the previous subsection, the claim   is  proved if we  show the existence of a pair of functions $\k$ and $\t$ of class $\cC^1$ satisfying \eqref{system} on  $(0, s_o]$. To prove this,
consider the equation and the initial conditions
\beq \label{5.18} \vartheta'(t) = \frac{\tan(\vartheta({\vartheta(t)}))}
{\sin(\vartheta(t) )}\ ,\qquad \lim_{t\to 0^+} \vartheta(t) = 0\ ,\ \lim_{t\to 0^+} \vartheta'(t) = A\eeq
on a $\cC^1$-function $\vartheta:(0,t_o] \to \bR$ and notice that if $\vartheta$ is a  solution to \eqref{5.18} with  $ \vartheta'(t) > 0$ at all points, then $\left(\t_s = R \ \vartheta^{-1}\left(\frac{s}{R}\right), \k_s = \frac{1}{R} \cot\left(\vartheta\left(\frac{s}{R}\right)\right)\right)$ is a $\cC^1$-solution of \eqref{system} on $(0, s_o = R\, \vartheta(t_o)]$ satisfying \eqref{5.17} with $a = \frac{1}{A}$. Hence  the lemma is proved if we  show that there exists a solution of \eqref{5.18} for any $0< A< 1$ and $t_o$ sufficiently small.\par
Given $0 < A < 1$, consider the following sequence of functions, defined by recursions on the same interval $[0, t_o] \subset \bR$ for some  $t_o$ to be determined later
\beq \vartheta_0(t) \= A t\ ,\qquad \vartheta_n(t) \= \int_0^t \frac{\tan\left( \vartheta_{n-1}(\vartheta_{n-1}(u))\right)}{\sin\left(\vartheta_{n-1}(u)\right)} d u\ .\eeq
We want to show that it is possible to determine $t_o$ such that the functions $\vartheta_n$  converge uniformly on $[0,t_o]$ to a solution of \eqref{5.18}. 
This is a direct consequence of  the following  claims. \par
\smallskip
\noindent {\bf Claim 1.}
{\it There exist 
$t_o \in (0,1)$ and $K < 1$ such that $\vartheta_n(t)$ is defined for any $n$ and any $t \in [0,t_o]$  and satisfies }
\beq  \label{cris} 
0 < \vartheta_n(t) \leq K t\, \qquad \text{for any}Ê\ t > 0\ .\eeq
\par
\smallskip
\noindent {\bf Claim 2.}   {\it $\lim_{t \to 0^+} \vartheta'_n(t) = A$ for any $n$. }\par
\noindent {\bf Claim 3.}   {\it $\vartheta''_n(t) \geq 0$ for any $n$ and any $t \in (0,t_o)$. }\par
\smallskip
\noindent {\bf Claim 4.} {\it $\vartheta_{n-1}(t) \leq  \vartheta_n(t)$ for any $n \geq 1$ and any $t \in [0,t_o]$.}\par
\smallskip
\noindent {\bf Claim 5.} {\it The sequence $\{\vartheta_n\}$ converges uniformly to a function $\vartheta:[0, t_o] \to \bR$ which is a solution of \eqref{5.18}}.\par
\par
\smallskip
Let us now proceed with the proofs of such claims.\par
\medskip
\noindent{\it Proof of Claim 1.} Consider a value $ B> 0$  such that 
\beq \label{cris1} A^5 + \frac{A^3}{2} + A^2 B - \frac{3 B}{2} < 0\eeq
and notice that the function 
$$F(x) = \tan\left( A x + \frac{B}{3} x^3\right) - \left( A + \frac{B}{2A^2} x^2\right) \sin\left( x\right)$$
is such that $F(0) = F'(0) = F''(0) = 0$ and 
\beq   F'''(0) = \frac{2}{A^2 } \left(A^5 + \frac{A^3}{2} + A^2 B - \frac{3  B}{2}\right) < 0\eeq
Therefore there exists an open interval $I = (0, \varepsilon)$,  on which $F(x) <  0$ and consequently
\beq \frac{ \tan\left( A x + \frac{B}{3 } x^3\right)}{ \sin\left( x\right)}  < \left( A + \frac{B}{2A^2} x^2\right)\ . \eeq
Now, let $t_o > 0$ satisfy: 
\begin{itemize}
\item[i)] $t_o \leq  \varepsilon$; 
\item[ii)] $\frac{1}{2} + \frac{B}{3A} t^2_o + \frac{B^2}{18 A^2} t_o^4 < 1$; 
\item[iii)] $A + B t_o^2 < 1$.
\end{itemize}
We now  show by induction  that for any $n$ and $t \in (0, t_o]$
\beq \label{crisbis} 0 <Ê\vartheta_n(t) \leq A t + \frac{B}{3} t^3 \ .\eeq
First of all, for $n = 0$, we have  that for any $t > 0$
\beq 0 < \vartheta_0(t) = A t \leq A t + B \frac{t^3}{3} \ .\eeq
Secondly, assume that $\vartheta_{n-1}(t)$ is defined for any $t \in [0, t_o]$ and that  $0 < \vartheta_{n-1}(t) \leq A t + B \frac{t^3}{3}$ for any $t \in (0, t_o]$. Then,  by (iii), we have that $\vartheta_{n-1}(t)   \leq A t + \frac{B t^3}{3} \leq t$ in $[0, t_o]$ and hence that $\vartheta_{n-1}([0,t_o]) \subset [0, t_o]$. It follows that 
the function $\vartheta_n$ is well defined in $[0, t_o]$ and $\cC^1$ with $ \vartheta'_t > 0$ for all $t> 0$.\par
From (i) and (ii) and inductive hypothesis, it follows that for any $t \in (0,t_o)$
$$ \vartheta'_n(t) = \frac{\tan\left(\vartheta_{n-1}(\vartheta_{n-1}(t))\right)}{\sin\left(\vartheta_{n-1}(t)\right)} \leq
\frac{\tan\left(A \vartheta_{n-1}(t) + \frac{B}{3}\vartheta_{n-1}^3(t)\right)}{\sin\left(\vartheta_{n-1}(t)\right)} \leq
$$
$$ \leq A + \frac{B}{2 A^2} \vartheta_{n-1}^2(t) \leq A + \frac{B}{2 A^2} \left( A^2 t^2 + \frac{2}{3} A B t^4 + B^2 \frac{t^6}{9}\right) \leq $$
\beq \label{ultima} \leq A + B t^2 \left( \frac{1}{2} + \frac{1}{3}ÿ\frac{B t^2}{A}  + \frac{B^2 t^4}{18 A^2}\right) \leq A + B t^2 \eeq
 From this,  the inequality \eqref{crisbis} follows by integration. Now,  by \eqref{crisbis} and the assumption (iii) on $t_o$ and  setting  $K = A + \frac{B t^2_o}{3} < 1$, we  obtain the inequality \eqref{cris}, namely 
$\vartheta_n(t)  \leq t \left(A + \frac{B t_o^2}{3}\right) = K t $.\par
\bigskip
\noindent{\it Proof of Claim 2.}  By construction,  $\vartheta'_0(0) = A$ and if $\lim_{t \to 0^+} \vartheta_{n-1}'(t) = A$, then 
$$\lim_{t \to 0^+} \vartheta_n'(t) = \lim_{t \to 0^+} \frac{\tan\left(\vartheta_{n-1}(\vartheta_{n-1}(t))\right)}{\sin\left(\vartheta_{n-1}(t)\right)}=  \lim_{t \to 0^+} \frac{\vartheta_{n-1}(\vartheta_{n-1}(t))}{\vartheta_{n-1}(t)} = $$
$$ = \lim_{t \to 0^+} \frac{\vartheta_{n-1}(\vartheta_{n-1}(t)) - \vartheta_{n-1}(0)}{\vartheta_{n-1}(t)} = \lim_{\t \to 0^+}  
\frac{\vartheta_{n-1}(\t) - \vartheta_{n-1}(0)}{\t} = $$
$$ = \lim_{\t \to 0^+} \vartheta'_{n-1}(\t) = A\ .$$
The claim follows by induction on $n$.\par
\bigskip
\noindent{\it Proof of Claim 3.} Also this claim is proved by induction. When $n = 0$, the claim is trivial since  $\vartheta_0''(t) = 0$.  Assume now that $\vartheta''_{n-1} \geq 0$ on $(0, t_o)$ and observe that 
$$\vartheta''_n(t) = \frac{d}{dt}  \frac{\tan\left(\vartheta_{n-1}(\vartheta_{n-1}(t))\right)}{\sin\left(\vartheta_{n-1}(t)\right)} = $$
$$ = \frac{  \frac{1}{\cos^2\left(\vartheta_{n-1}(\vartheta_{n-1}(t))\right)} \vartheta'_{n-1}(\vartheta_{n-1}(t)) \vartheta'_{n-1}(t)  }{ \sin\left(\vartheta_{n-1}(t)\right)   } -$$
$$ -  \frac{ \cos\left(\vartheta_{n-1}(t)\right) \tan\left(\vartheta_{n-1}(\vartheta_{n-1}(t))\right) \vartheta'_{n-1}(t)   }{  \sin^2\left(\vartheta_{n-1}(t)\right)} = $$
\beq \label{5.26} =  \frac{\vartheta'_{n-1}(t) \cos\left(\vartheta_{n-1}(t)\right)}{ \cos^2\left(\vartheta_{n-1}(\vartheta_{n-1}(t))\right) \sin^2\left(\vartheta_{n-1}(t)\right)}\cdot \cF^{(n-1)}(\vartheta_{n-1}(s)) \eeq
where, for any $m$,  we denote 
 \beq \label{fn2} \cF^{(m)}(s) \= \vartheta'_{m}(s) \tan\left(s\right) - \sin\left( \vartheta_{m}(s)\right) \cos\left( \vartheta_{m}(s)\right)\ . \eeq
 Using the assumption  $\vartheta''_{n-1}(s) \geq 0$ for $s \in (0,t_o]$ and the fact that  $\vartheta'_{n-1}(s) > 0$, we have that
$$ \cF^{(n-1)}{}'(s) = \vartheta''_{n-1}(s) \tan(s) + \vartheta'_{n-1}(s) \frac{1}{\cos^2(s)} - $$
$$ - \vartheta'_{n-1}(s) (\cos^2(\vartheta_{n-1}(s)) - \sin^2(\vartheta_{n-1}(s))) = $$
$$ = \vartheta''_{n-1}(s) \tan(s) + \frac{\vartheta'_{n-1}(s)}{\cos^2(s)} \left( 1 - \cos^2(s) \cos(2 \,\vartheta_{n-1}(s))\right)
\geq 0\ .$$
Since $\cF^{(m)}(0) = 0$ for any $m$, it follows that $\cF^{(n-1)}(s) \geq 0$ for any $s \in (0, t_o]$ and  hence that $\cF^{(n-1)}(\vartheta_{n-1}(t)) \geq 0$ for any $t \in (0, t_o]$. From this and \eqref{5.26}, we 
get $\vartheta''_n(t) \geq 0$ as needed.\par 
\bigskip
\noindent{\it Proof of Claim 4.}  First of all, notice that for any $t \in (0,t_o]$
$$\vartheta'_1(t) = \frac{\tan(\vartheta_0(\vartheta_0(t)))}{\sin(\vartheta_0(t))} = \frac{\tan(A ^2 t)}{\sin(A t)} \geq \frac{A^2 t}{A t}  = A =  \vartheta_0'(t)\ .$$
Hence, by integration,  $\vartheta_1(t) \geq \vartheta_0(t)$ for any $t \in (0,t_o]$. Let us  assume that  $\vartheta_{n-2} \leq \vartheta_{n-1}$ at all points of $[0,t_o]$. In order to prove the claim by an inductive argument, we only need to check that $\vartheta'_{n-1} \leq \vartheta'_n$ on the same interval $[0, t_o]$. To see this, we notice that 
$$\vartheta'_n(t) - \vartheta'_{n-1}(t) = \frac{\tan(\vartheta_{n-1}(\vartheta_{n-1}(t)))}{\sin (\vartheta_{n-1}(t))} - 
\frac{\tan(\vartheta_{n-2}(\vartheta_{n-2}(t)))}{\sin (\vartheta_{n-2}(t))} \overset{\vartheta_{n-1} \geq \vartheta_{n-2}} \geq 
$$
$$ \geq \frac{\tan(\vartheta_{n-2}(\vartheta_{n-1}(t))}{\sin(\vartheta_{n-1}(t))} - 
 \frac{\tan(\vartheta_{n-2}(\vartheta_{n-2}(t))}{\sin(\vartheta_{n-2}(t))} =$$
 $$ = \frac{d}{d s} \left.\frac{\tan(\vartheta_{n-2}(s)}{\sin(s)}\right|_{s = \wt s} \left(\vartheta_{n-1}(t) - \vartheta_{n-2}(t)\right)$$
 for some $\wt s\in (\vartheta_{n-2}(t), \vartheta_{n-1}(t))$. On the other hand 
 $$\frac{d}{d s} \frac{\tan(\vartheta_{n-2}(s))}{\sin(s)} =$$
 $$ =  \frac{\cos s}{\cos^2(\vartheta_{n-2}(s)) \sin^2 s} \left\{ \vartheta'_{n-2}(s) \tan(s) - \sin( \vartheta_{n-2}(s)) \cos( \vartheta_{n-2}(s))\right\} =$$
 $$ =  \frac{\cos s}{\cos^2(\vartheta_{n-2}(s)) \sin^2 s}  \cF^{(n-2)}(s)$$
where $\cF^{(n-2)}(s)$ is  as  defined in \eqref{fn2}.  In the proof of the previous claim,  we showed that $\cF^{(m)} \geq 0$ for any $m \geq 1$ and from this  we conclude that $\vartheta'_{n-1} \leq \vartheta'_n$ as needed. \par
\par
\bigskip
\noindent{\it Proof of Claim 5.}
For the proof of this claim, we need the following properties which are consequences  of the  previous claims:
\begin{itemize}
\item[i)] $\vartheta'_{n-1}(\vartheta_{n-1}(t)) \leq \vartheta_{n-1}'(t)$; 
\item[ii)] $\cos^2(\vartheta_{n-1}(\vartheta_{n-1}(t))) \geq \cos^2(\vartheta_{n-1}(t))$; 
\item[iii)] $\vartheta'_{n-1}(t) \leq 1$; 
\item[iv)] $\frac{1}{\cos^2(\vartheta^2_{n-1}(t))} \leq \frac{1}{\cos^2(\vartheta_{n-1}^2(t_o))} = K$.
\end{itemize} 
Using these relations, for any $t \leq t_o$ we have that  (here ``$\vartheta_m^2(t)$'' stands for ``$\vartheta_m(\vartheta_m(t))$''): 
$$\vartheta''_n(t)  = \frac{\vartheta'_{n-1}(t) \cos(\vartheta_{n-1}(t))}{\cos^2(\vartheta^2_{n-1}(t) )\sin^2(\vartheta_{n-1}(t))} \cdot$$
$$\cdot \left\{ \vartheta'_{n-1}(\vartheta_{n-1}(t)) \tan(\vartheta_{n-1}(t)) - 
\sin(\vartheta^2_{n-1}(t)) \cos (\vartheta^2_{n-1}(t)) \right\} \leq $$
$$\leq  \frac{K}{\sin^2 \vartheta_{n-1}(t)} \left\{ \vartheta'_{n-1} (\vartheta_{n-1}(t))\tan(\vartheta_{n-1}(t)) - 
\sin(\vartheta^2_{n-1}(t)) \cos(\vartheta^2_{n-1}(t) )\right\} = $$
$$ = K \left\{ \vartheta'_{n-1}(\vartheta_{n-1}(t)) \frac{1}{\sin(\vartheta_{n-1}(t)) \cos(\vartheta_{n-1}(t))} - 
\frac{\vartheta'_n(t)  \cos^2(\vartheta_{n-1}^2(t))}{\sin(\vartheta_{n-1}(t))}\right\} =$$
$${= \frac{K \left\{ \vartheta'_{n-1}(\vartheta_{n-1}(t)) - \vartheta_n' (t)\cos (\vartheta_{n-1}(t)) \cos^2 (\vartheta^2_{n-1}(t))\right\} }{\sin (\vartheta_{n-1}(t) )\cos (\vartheta_{n-1}(t))}}\leq
 $$
$$\overset{(ii)} \leq \frac{K'}{\sin (\vartheta_{n-1}(t))} \left\{  \vartheta'_{n-1}(\vartheta_{n-1}(t)) - \vartheta'_n(t) \cos^3 (\vartheta_{n-1}(t))\right\}    \leq  $$
$$\overset{\vartheta'_{n-1}(t) \leq \vartheta'_n(t)} \leq \frac{K' \vartheta'_n(t)}{\sin( \vartheta_{n-1}(t))} \left\{ 1 - \cos^3 (\vartheta_{n-1}(t))\right\} \leq $$
\beq \label{uffaa} \leq 3 K'  \frac{ 1 - \cos (\vartheta_{n-1}(t))}{\sin (\vartheta_{n-1}(t))} =
3 K'  \tan\left(\frac{\vartheta_{n-1}(t)}{2}\right) \leq 3 K' \tan\left(\frac{t_o}{2}\right)
 \ ,\eeq
 i.e.  the functions $\vartheta''_n|_{[0, t_o]}$ are uniformly bounded. From this we get that  the  $\vartheta'_n|_{[0, t_o]}$ are uniformly bounded and equicontinuous. Since 
 the sequence $\{ \vartheta'_n|_{[0, t_o]}\}$ is monotone,  the sequence $\{ \vartheta_n|_{[0, t_o]}\}$ uniformly converges  to a $\cC^1$-function $\vartheta: [0, t_o] \to \bR$  satisfying 
$$\lim_{n \to \infty}Ê\vartheta'_n(t) = \vartheta'(t)\qquad \text{for any} \ t \in [0, t_o]\ .$$
We leave to the reader the simple task of checking that  the sequence $\vartheta_n(\vartheta_n(t))$ uniformly converges to $\vartheta(\vartheta(t))$, from which follows  that the limit function $\vartheta(t)$
is indeed a solution to the differential problem \eqref{5.18}.
\end{pf}
\bigskip
\subsection{Existence of   self-involutes on the spheres} \hfill\par
\begin{theo} \label{theoremonselfinvolutes} Let   $\h: [0, L] \to S^2_R$ be an almost self-involute on $S^2_R$ of class $\cC^3$, parameterized by arc length,  such that 
\beq \label{5.17bis} \lim_{s\to 0^+} \t_s = 0\ ,\quad  \lim_{s\to 0^+} \dot \t_s = a\ ,\quad  \lim_{s\to 0^+} \k_s = + \infty\ . \eeq
Then $\h$ is  self-involute if and only if $a$ is the unique solution of the equation 
\beq \label{fundamentalequation} a = Êe^{\frac{3 \pi}{2a}} \ .\eeq
\end{theo}
\begin{pf} Given an almost self-involute $\h: [0,L] \to S^2_R$,  let 
$\wt \h$ be the initial arc of the  involute that is congruent to $\h$. With no loss of generality, we assume that $\h_0 = (0,0,R)$. By definition, we have that $\wt \h_0 = \h_0 = (0, 0, R)$ and
that there exists an orthogonal matrix 
\beq \label{rotation}ÊA^{(\h)} = \left(\begin{matrix}  \cos \phi & - \sin \phi\ & 0\\   \sin \phi & \cos \phi & 0 \\
0 & 0 & 1\end{matrix} \right)\qquad \text{or}\qquad A^{(\h)} = \left(\begin{matrix}  \cos \phi & \sin \phi\ & 0\\   \sin \phi & - \cos \phi & 0 \\
0 & 0 & 1\end{matrix} \right)\eeq
such that $ \wt \h_s = A^{(\h)} \cdot \h_s$ for any $s \in [0,L]$. \par
We also identify   $E^2$ with the plane $\{\ x^3 = 1\ \}$ in $\bR^3$ and also in this case 
 for any 
almost self-involute $\h: [0,L] \to E^2 = \{\ x^3 = 1 \ \}$ with $\h_0 = (0,0,1)$, we denote by $\wt \h$ its involute and by $A^{(\h)}$ the  orthogonal matrix as in \eqref{rotation} such that  $ \wt \h_s = A^{(\h)} \cdot \h_s$. \par
Finally, for any almost self-involute $\h$ on $S^2_R$ or $E^2$, we 
call {\it fundamental   pair of $\h$\/} the couple $(A, a)$ given by  $A = A^{(\h)}$ and   $ a = \lim_{s\to 0^+} \dot \t_s$, where $\t$ is defined in  \eqref{system} or $\eqref{system1}$, respectively. Clearly,  $\h$ is self-involute if and only if $A = I$.\par
\medskip
The proof is a consequence of the existence of a canonical correspondence between  any almost self-involute $\h$ on  $S^2_R$ with fundamental  pair $(a, A)$ and  an almost self-involute $\h^{(\infty)}_s$ on $E^2$, determined up to Euclidean isometries, having  the same fundamental pair of $\h$ and congruent to the curve
 \small 
\beq\label{MP}\h^{(\infty)}_s = \left(\frac{s}{\sqrt{1 + a^2}} \cos\left( a \log\left(\frac{s}{\sqrt{1 + a^2}}\right) \right), 
- \frac{s}{\sqrt{1 + a^2}} \sin\left( a \log\left(\frac{s}{\sqrt{1 + a^2}}\right) \right), 1\right) \eeq 
\normalsize
which is  a parameterization by arc length of the curve 
$$\g_t = (e^{- \frac{t}{a}} \cos t, e^{- \frac{t}{a}} \sin t, 1)\ ,$$
studied   by  Manselli and Pucci in \cite{MP1}. In that paper, it is proved that  {\it $\h^{(\infty)}$ is self-involute\/}Ê (i.e. with a fundamental pair  $(I, a)$) {\it if and only if $a$ is solution of \eqref{fundamentalequation}}.  Since $\h^{(\infty)}$ and  $\h$ have the same fundamental pair,  the conclusion  follows.\par
\medskip
Let us prove the existence of   the correspondence  $\h \longmapsto \h^{(\infty)}$ described above. For any $0 \neq \l \in \bR$,  let 
 \beq\label{dilated} \h^{(\l)}: [0, L] \to S^2_{\l R}\ ,\qquad  \h^{(\l)}_s \= \l\ \h \left(\frac{s}{\l}\right)\ , \eeq
 which is the initial  arc of length $L$ (parameterized by arc length) of the dilatation  of $\h$ by $\l$.  
From definitions, one can check that $\h^{(\l)}$ is  an almost self-involute of $S^2_{\l R}$, with  
involute $\wt{\h^{(\l)}}Ê= \wt \h^{(\l)} = A \cdot \h^{(\l)}$ and with functions $\k^{(\l)}_s$ and $\t^{(\l)}_s$ given by 
\beq \k^{(\l)}_s = \frac{1}{\l} \k_{\frac{s}{\l}} \ , \qquad  \t^{(\l)}_s = \l \t_{\frac{s}{\l}} \ . \eeq
In particular, $\lim_{s \to 0} \dot \t^{(\l)}_s = \lim_{s \to 0} \dot \t_s = a$ and the fundamental pair $(A, a)$ is the same for all curves $\h^{(\l)}$. 
For all $\l$ sufficiently large,    $\h^{(\l)}: [0, L] \to S^2_{\l R}$ is  included in the upper hemisphere
and  it can be identified with its image in  $ \{ x^3 = \l R\}$ by  the projection $\pi^{(\l)}:  S^2_{\l R(+)} \to  \{ x^3 = \l R\} \simeq E^2$
of center the origin.  
\par
\smallskip
The correspondence   we are looking for  is based on the  following lemma.
\par
\smallskip
\begin{lem} There exists a sequence $\l_n \to + \infty$ such that the curves $\h^{(\l_n)}$ converge on $(0, L]$, uniformly on compacta together with their first and second derivatives,   to a $\cC^2$-curve $\h^{(\infty)}$. 
\end{lem}
\begin{pf} To prove this, first of all notice that, being of length $L$ on the sphere and obtained via  the projection $\pi^{(\l)}$,  any curve $\h^{(\l)}$  starts from $x_o = (0,0, \l R) (\simeq (0,0,1) \in E^2)$ and it  is contained in the closed disk $\overline{D_r(x_o)}$ of radius $r = \l R \tan\left(\frac{L}{\l R}\right) < 2 L$ for all $\l$ sufficiently large. Secondly,  let us denote by $g^{(\l)}$ the metric on $E^2 \simeq \{\ x^3 = \l R\ \}$ defined as   push-forward by the projection $\pi^{(\l)}$ of the metric of $S^2_{\l R}$. Notice that   on any closed disc $\overline{D_{r_o}(x_o)}$, the metric $g^{(\l)}$ converges uniformly to the Euclidean metric $g_o$  together with  all  derivatives. Now, 
by construction, for any $\l$ and $s \in [0, L]$, we have that 
$g^{(\l)}(\dot \h^{(\l)}_s, \dot \h^{(\l)}_s) = 1$
and hence $|\dot \h^{(\l)}_s| = \sqrt{ g_o\left(\dot \h^{(\l)}_s, \dot \h^{(\l)}_s\right)}$ is uniformly bounded for all $\l$ sufficiently large. A similar argument shows that also  the normal vectors
$\Norm^{(\l)}_s$ of the curve $\h^{(\l)}$ (orthogonal to the $\Tang^{(\l)}_s = \dot \h_s$ w.r.t. $g^{(\l)}$) are uniformly bounded.
\medskip
On the other hand,   from \eqref{system},  the fact that $\t$ is monotone and that  $\lim_{s \to 0^+} \t_s = 0$,  one has that for any fixed $0 < \varepsilon_o < L$ and any $s \in [\varepsilon_o, L]$
 \beq\label{limitcurvature}Ê \lim_{\l \to +\infty} \k^{(\l)}_{s} = \lim_{\l \to  +\infty} \frac{1}{\l} \k_{\frac{s}{\l}} = 
  \lim_{\l \to + \infty} \frac{1}{\l R} \cot\left(\frac{\t^{-1}\left(\frac{s}{\l}\right)}{R}\right) =$$
  $$ =   \lim_{\mu \to  0^+} \frac{\mu}{\t^{-1}\left(\mu s\right)}\overset{\text{l'H\^opital}}= \frac{a}{s} \leq  \frac{a}{\varepsilon_o}\eeq
  and with similar computations
  $$\lim_{\l \to +\infty} \dot \k^{(\l)}_{s} = -  \frac{a}{s^2} \geq -  \frac{a}{\varepsilon_o^2}\ .$$
  From this and \eqref{Frenetformulae},  it follows that the covariant derivatives $\n_{\dot \h^{(\l)}} \Tang^{(\l)}_s$ and $\n^2_{\dot \h^{(\l)}} \Tang^{(\l)}_s$  are uniformly bounded in any given interval $[\varepsilon_o, L]$. Considering  the explicit expression of such  covariant derivatives in terms of the Christoffel symbols of $g^{(\l)}$ and of the derivatives $\ddot \h^{(\l)}_s$  and $\dddot \h^{(\l)}_s$, one can directly check that on any given interval $ [\varepsilon_o, L]$, the curves $\h^{(\l)}$ are uniformly bounded in $\cC^3$-norm.  From this, the lemma follows.
  \end{pf}
  \medskip
 Using definitions and convergence in  $\cC^2$, one can check that if $(A, a)$ is the fundamental pair of $\h: [0, L] \to S^2_R$, then 
 the limit curve $\h^{(\infty)}$   and  $\wt \h^{(\infty)} = A \cdot \h^{(\infty)}$ satisfy the relation  \eqref{involute}, i.e. $\h^{(\infty)}$ is almost self-involute. Moreover, by \eqref{limitcurvature}, we see that  the curvature  of $\h^{(\infty)}$ is given by $\k^{(\infty)}_s = \frac{a}{s}$ and the associated function  is $\t^{(\infty)}_s = a s$.  In particular,   the fundamental pair of $\h^{(\infty)}$ is $(A, a)$,  the same of $\h$. Notice also that  since  the curve \eqref{MP} has curvature function given by  $\k_s = \frac{ a}{s}$, by the  Fundamental Theorem of Plane Curves,  $\h^{(\infty)}$ is congruent to \eqref{MP} and this concludes the proof that   the  correspondence  $\h \longmapsto \h^{(\infty)}$ has  all stated properties.
\end{pf}
From this and Theorem \ref{existenceselfinvolutes}, the next corollary follows  immediately. 
\begin{cor} \label{thecorollary}ÊThere exists  $L_o > 0$ such that for any $L < L_o$  there exists a  self-involute curve  of length $L$ on the unit sphere $S^2$. 
\end{cor}
\bigskip
\section{$\cG$-curves on  spheres  with the ``maximal length property'' }
\setcounter{equation}{0}
\medskip
This section is devoted to show the existence  of
$\cG$-curves $\g:[0, L] \to \cS = S^2_+$ realizing  the ``maximal length property'', i.e. such that for any $s \in [0,L]$, 
the length $\ell_\cS(\g|_{[0,s]}) = s$  coincides with   the largest  possible value according to Theorem \ref{maintheorem1}, namely  $\ell(\g|_{[0,s]}) = Ê\gp(s) $. \par
\smallskip
\begin{theo}\label{finaltheorem}  If $\h:[0,L] \to \cS = S^2_+$   is a self-involute, then it is a $\cG$-curve with the maximal length property, i. e.  such that 
\beq \label{conditionbis}Ê\gp(s) =  \ell_\cS(\g|_{[0,s]}) \ , \qquad \text{for any} \ s \in [0,L]\ .\eeq
\end{theo}
\begin{pf} In the following, we constantly identify the self-involute $\h$ with its image in  $ \{ x^3 = 1\} \simeq E^2$, determined  by  the projection $\pi:  S^2_{(+)} \to  \{ x^3 = 1\} $ of center the origin.  We also denote by $g$ the metric  on $E^2$ defined as   push-forward by the projection $\pi$ of the metric of $S^2$ and by $g_o = d x^1 \otimes dx^1 + dx^2 \otimes dx^2$ the standard Euclidean metric. 
\par
\smallskip
A simple computation shows
$$g = g_{ij}Êdx^i \otimes dx^j\ \ \text{with}\ \ 
\left(\smallmatrix g_{11} & g_{12}\\ \phantom{a}\\ g_{21} & g_{22} \endsmallmatrixÊ\right)
= \smallmatrix \frac{1}{(1 +  (x^1)^2 + (x^2)^2)^2}\endsmallmatrix \left( \smallmatrix 
1 + (x^2)^2  & -  x^1 x^2 \\
 - x^1 x^2 & 
1 + (x^1)^2\endsmallmatrix \right)\ \  .$$
From this,  by well-known formulae, we may 
 express in term of $\h_s = (\h_s^1, \h_s^2)$ and its derivatives the following objects:
\begin{itemize}
\item[--] the Christoffel symbols $\left.\G_{ij}^{\ k}\right|_{\h_s} =\left. \frac{1}{2} g^{km}\left(\frac{\partial g_{mj}}{\partial x^i} + \frac{\partial g_{im}}{\partial x^j} - \frac{\partial g_{ij}}{\partial x^m} \right)\right|_{\h_s}$; 
\item[--] the covariant derivatives 
$\nabla_{\Tang_s} \Tang_s = \left(\ddot \h_s^i +\left. \G_{jk}^{\ i} \right|_{\h_s} \dot \h^j_s \dot \h^k_s \right)\frac{\partial}{\partial x^i}$; 
\item[--]the geodesic curvature $k_s = \sqrt{g(\nabla_{\Tang_s} \Tang_s , \nabla_{\Tang_s} \Tang_s )} = g(\nabla_{\Tang_s} \Tang_s, \Norm_s)$.
\end{itemize} 
Since $\k_s > 0$, also the Euclidean curvature $\k^E_s$ of $\h_s$ is positive. This can be checked as follows. Recall that $\ddot \h_s = g_o(\dot \h_s, \dot \h_s) \k^E_s \Norm^E_s + \l_s \dot \h_s$ for some function $\l_s$ and with $\Norm^E_s$ denoting the  unit normal vector in the Euclidean sense. Since $g(\dot \h_s, \Norm_s) = g(\Tang_s, \Norm_s) = 0$ and using the expression for $ \nabla_{\Tang_s} \Tang_s $,  
$$ g_o(\dot \h_s, \dot \h_s) \k^E_s g(\Norm^E_s, \Norm_s) = g(\ddot \h_s, \Norm_s) =  \k_s  - g\left(\left. \G_{jk}^{\ i} \right|_{\h_s} \dot \h^j_s \dot \h^k_s\frac{\partial}{\partial x^i}, \Norm_s \right). $$
Using properties of projectively flat connections or just by a direct computation, one can see that the vector $v_s = \left. \G_{jk}^{\ i} \right|_{\h_s} \dot \h^j_s \dot \h^k_s \frac{\partial}{\partial x^i}$ is proportional to $\dot \h_s = \Tang_s$. From this,  it follows that 
$\k^E_s  =  \frac{\k_s}{g_o(\dot \h_s, \dot \h_s)  g(\Norm^E_s, \Norm_s)}$. Since $\Norm^E_s$ and $\Norm_s$ lie on the same side w.r.t. $\Tang_s$, $g(\Norm^E_s, \Norm_s) > 0$ and $\k^E_s$ is positive.  \par
\smallskip
Being $\k^E > 0$,   the ``angle'' function (i.e. the Euclidean angle between $\dot \h_s$ and the $x^1$-axis), computable by 
$$\varphi: [0, L] \to \bR\ ,\  \varphi \=  \int_L^s \k^E_u d u + \varphi_L\ ,\quad \text{with}\ \varphi_L = \wh{\dot \h_L \frac{\partial}{\partial x^1}}  ,$$
is monotone increasing.  \par
\bigskip
Now, for any $s \in (0,L]$,  let us consider the following notation: 
\begin{itemize}
\item[--] $\ell_1^{(s)}$ denotes the tangent line to $\h$ at the point $x = \h_s$, while $\ell_2^{(s)}$ denotes 
the line through $x$ and parallel to the vector $\Norm_s$; recall that these lines, up to   reparameterizations, 
are    geodesics for both the Euclidean metric $g_o$  and  the spherical metric $g$;  
\item[--] $\g^{(s)}$ is the closed, piecewise $\cC^1$ curve,  formed by  the arc $\h|_{[\t^{-1}(s), s]}$ and the segment joining  $\h_s$ and $\h_{\t^{-1}(s)}$; 
\item[--]  $\rho^{(s)}$ is the total rotation of   $\g^{(s)}$, i.e.  the multiple of $2\pi$ defined by 
\beq \label{ns} \rho^{(s)} = \varphi_s - \varphi_{\t^{-1}(s)}  + \wh{ \Norm_s \Tang_s}\eeq
and $m_s = \rho^{(s)}/ 2 \pi$.   Being $\varphi$ monotone increasing, $\rho^{(s)}$ coincides with the total 
curvature of the curve  $\g^{(s)}$. From this,  by  Fenchel's theorem for piecewise differentiable curves (\cite{Fe}; see also \cite{Mi, Ae}), $\g^{(s)}$ is a simple convex curve if and only if $m_s = 1$. 
\end{itemize}
We claim that  for any $s$,  $\ell^{(s)}_1$ and $\ell^{(s)}_2$ are support lines for the arc $\h|_{[0,s]}$ and that  
$\g^{(s)}$ is 
a simple and convex curve. 
This immediately implies that $\h$ is a $\cG$-curve and that $\g^{(s)}$ is the boundary of the 
convex hull of  $\h|_{[0, s]}$. Being $\h_s$  self-involute, it follows that 
$$\gp(s) = \ell_{\cS}(\h|_{[\t^{-1}(s), s]}) +  \ell_{\cS}([\h_{\t^{-1}(s)}, \h_s])  = s - \t^{-1}(s) + \t^{-1}(s) = s$$
i.e.  $\h_s$ satisfies  \eqref{conditionbis} at any $s$.\par
\medskip
 The proof of the theorem is therefore  a direct consequence  of the following three claims.\par
 \medskip
 \noindent{\bf Claim 1: }  {\it  $m_s= 1$ for any $s$, i.e. any closed curve $\g^{(s)}$  is simple and convex}. \par
 \smallskip
 \noindent
 In fact, from \eqref{ns}, the map  $s \mapsto m_s$ is  continuous  and  therefore   constant.  Moreover,   for any $\l > 1$, the arc $\h|_{[0, \frac{L}{\l}]}$ is homothetic to the curve $\h^{(\l)}$ described in \eqref{dilated} (which we also identify with the corresponding projected curve on $E^2$) and, by the proof of 
Theorem \ref{theoremonselfinvolutes},  we may  choose $\l$ so large   that  $\h^{(\l)}$ is arbitrarily close
in $\cC^2$-norm to the self-involute $\h^{(\infty)}$. By the results in \cite{MP}, $\h^{(\infty)}$ is a  $\cG$-curve  of $E^2$. So, 
if we denote by $\g^{(s|\infty)}$ the closed curved formed by the segment joining $\h^{(\infty)}_s$ and $\h^{(\infty)}_{(\t^{(\infty)})^{-1}(s)}$ and the arc  $\h^{(\infty)}|_{[(\t^{(\infty)})^{-1}(s), s]}$, it is immediate  to realize that $\g^{(s| \infty)}$ is simple and convex, that is $m^{(\infty)}_s = 1$. Since 
 $\g^{(s)}$ is homothetic to the piecewise $\cC^2$ closed curve $\g^{(s|\l)}$,  close to  $\g^{(s|\infty)}$ in $\cC^2$-norm, it follows that also $m_s = 1$ for any $0 < s \leq \frac{L}{\l}$ and hence for all values of $s$. 
\par
\medskip
\noindent{\bf Claim 2:} {\it  $\ell_1^{(s)}$ is a support line for $\h|_{[0,s]}$}. 
\par
\smallskip
\noindent To see this, notice that, being $\g^{(s)}$ closed and convex,  $\ell_1^{(s)}$ is a support line for $\g^{(s)}$ and hence for 
$\h_{[\t^{-1}(s), s]}$. On the other hand, the spherical distance  $d_{\cS}(\t^{-1}(s), \ell_1^{(s)}) $ is equal 
to the length   of the segment joining $\h_{\t^{-1}(s)}$ and $\h_s$, because it  lies in a line  which is $g$-orthogonal to $\ell_1^{(s)}$. The length of this segment is equal to 
$\t^{-1}(s)$ by the definition of self-involute. This  length is also equal to the length  of  the arc $\h_{[0,\t^{-1}(s)]}$.  Hence,  this arc lies
in  the 
same half-plane of $\h_{[\t^{-1}(s), s]}$ and  $\ell^{(s)}_1$ is a support line for the  whole curve  $\h|_{[0,s]}$.\par
\medskip
\noindent{\bf Claim 3:} {\it  $\ell_2^{(s)}$ is a support line for $\h|_{[0,s]}$}. 
\par
\smallskip
\noindent As before, being $\g^{(s)}$ closed and convex,  $\ell_2^{(s)}$ is a support line for $\h_{[\t^{-1}(s), s]}$. On the other hand, by definition of self-involute,  $\ell_2^{(s)} = \ell_1^{(\t^{-1}(s))}$ and hence, by Claim 2,  it is also a support line for the arc $\h|_{[0,\t^{-1}(s)]}$. The two  arcs lie in the same half-plane,  because the Euclidean curvature of $\h$ is strictly positive at all points
 and the claim follows.  
\end{pf}

\end{document}